\documentclass{amsart}
\usepackage{ifthen}
\def\style{}
\ifthenelse{\equal{\style}{preprint}}{
\usepackage{dmo_preprint}
}{}
\usepackage{hyperref}
\usepackage{url}

\usepackage{amsmath,amssymb,amsthm,verbatim,paralist,graphicx}
\usepackage[active]{srcltx}
\usepackage[T1]{fontenc}
\usepackage{yfonts}
\usepackage{color}
\DeclareMathOperator{\rk}{r}
\DeclareMathOperator{\cl}{cl}
\DeclareMathOperator{\sep}{sep}
\DeclareMathOperator{\supp}{supp}
\DeclareMathOperator{\Dir}{Dir}

\newtheorem{theorem}{Theorem}
\newtheorem{prop}{Proposition}
\newtheorem{lemma}{Lemma}
\newtheorem{corollary}{Corollary}
\newtheorem{definition}{Definition}
\newtheorem{conjecture}{Conjecture}
\newtheorem{remark}[lemma]{Remark}

\begin{document}
\title{Mutations and (Non-)Euclideaness in oriented matroids} \author{Michael Wilhelmi}
\address{FernUniversit\"at in Hagen}
\email{\{Michael.Wilhelmi\}@FernUniversitaet-Hagen.de}

\ifthenelse{\equal{\style}{preprint}}
{ \DMOmathsubject{05B35, 05B25, 06C10, 51D20} \DMOkeywords{matroids, amalgams, embeddings, projective spaces} \DMOtitle{044.17}{Sticky matroids and
    Kantor's Conjecture} {Winfried
    Hochst\"attler\\Michael Wilhelmi}{\{Winfried.Hochstaettler, Michael.Wilhelmi\}@fernuni-hagen.de} }{}

\begin{abstract}
We call an oriented matroid {\em Mandel} if it has an extension in general position
which makes all programs with that extension Euclidean. If $L$ is the minimum number of mutations adjacent to an element of the groundset,
we call an oriented matroid {\em Las Vergnas} if  $L > 0$. 
If $\frak{O}_{property}$ is the class of oriented 
matroids having a certain property, it holds
\[ \frak{O} \supset \frak{O}_{Las Vergnas} \supset \frak{O}_{Mandel} \supset  \frak{O}_{Euclidean}  \supset  \frak{O}_{realizable}.\]
All these inclusions are proper, we give explicit proofs/examples for the parts of this chain that were not known. 

For realizable hyperplane arrangements of rank $\rk$ we have $L = \rk$ which was proved by Shannon.
Under the assumption that a (modified) intersection property holds we give an analogon to Shannons proof and show that 
uniform rank $4$ Euclidean oriented matroids with that property have $L = 4$. 
Using the fact that the lexicographic extension creates and destroys certain mutations, we show that
for Euclidean oriented matroids holds $L \ge 3$. 

We give a survey of preservation of Euclideaness and prove that Euclideaness remains after a certain type of mutation-flips.
This yields that a path in the mutation graph from a Euclidean oriented matroid to a totally non-Euclidean oriented matroid 
(which has no Euclidean oriented matroid programs) must have at least three mutation-flips. 
Finally, a minimal non-Euclidean or rank $4$ uniform oriented matroid is Mandel if it is connected to a Euclidean oriented matroid
via one mutation-flip, hence we get many examples for Non-Euclidean but Mandel oriented matroids and have $L \le 3$ for those of rank $4$. 
\end{abstract}

\maketitle

\section{introduction}

Arnaldo Mandel gave in his dissertation \cite{MandelDiss} the following conjecture.
\begin{conjecture}
Every oriented matroid $\mathcal{O}$ has an extension $g$ in general position (not a coloop)
such that $(\mathcal{O} \cup g, g, f)$ is a Euclidean oriented matroid program for all elements $f$ of $\mathcal{O}$ that are not loops or coloops.
\end{conjecture}
We call oriented matroids fulfilling that conjecture {\em Mandel}.
It was already known by Mandel, proved in his thesis, that in a Mandel oriented matroid each element must have
at least one adjacent simplicial tope / mutation. We call oriented matroids with that property {\em Las Vergnas}.
That fact is again shown in \cite{Knauer} and we show it also in this paper using more recent
language and avoiding the theory of partial cubes. 
Years after Mandels thesis three oriented matroids were found that are not Las Vergnas (the $R(20)$ by Richter Gebert given in \cite{Richter-Gebert}, Theorem 2.3
with $21$ elements,
 an oriented matroid similar constructed with $17$ elements by Bokoswki/Rohlfs in \cite{BokowskiMutationProblem}, Chapter 6, 
 and a $13$ element rank-$4$ oriented matroid by Tracy Hall in \cite{hall2004counterexamples}, Chapter 2.4 ). Hence Mandels conjecture turned out to be wrong. We will show in this paper that the following series of proper inclusions holds:
 \begin{equation}\label{equ:inclusions}
 \frak{O} \overset{\text{(1)}}{\supsetneq} \frak{O}_{Las Vergnas} \overset{\text{(2)}}{\supsetneq} \frak{O}_{Mandel} \overset{\text{(3)}}{\supsetneq}  \frak{O}_{Euclidean}  \overset{\text{(4)}}{\supsetneq}  \frak{O}_{realizable}
 \end{equation}
where $\frak{O}_{property}$ is the class of oriented 
matroids having a certain property. We give explicit proofs and examples for the properness where there are missing
(or assumed to be obvious but not given).

Inclusion (4) is classic (for the properness recall that all non-linear rank $3$ oriented matroids are Euclidean).
The three examples of mutation-free oriented matroids mentioned above show the properness of Inclusion (1), 
we considered already the $\supseteq$-part of Inclusion (2).
Inclusion (3) seems obvious and is mentioned in \cite{Knauer}. We show the $\supseteq$-part
by using an arbitrary lexicographic extension in general position in a Euclidean oriented matroid.
It serves as the desired extension because we showed in \cite{LexExt} that
lexicographic extensions preserve Euclideaness. We use this fact a lot in this paper and
 develop the theory of these extensions further connecting it with the theory of mutations.

Because a lexicographic extension yields an inseparable pair, it creates new mutations (between both elements of the pair), 
certain mutations will be destroyed (the new element cuts through the simplex) 
and others are preserved (where the pair is not involved). We use these results in several directions. 
First, we take the $R(20)$ which is non-Las Vergnas and has exactly one mutation-free element $f$. We extend that element lexicographically (which creates a new
mutation adjacent to $f$) and obtain a Las Vergnas oriented matroid. But this oriented matroid can not be Mandel
because the class $ \frak{O}_{Mandel}$ is minor closed (we show that here).
Hence, inclusion (2) is proper and the class $ \frak{O}_{LasVergnas}$ is not minor-closed.

Next, we consider the number $L$ of the minimum adjacent mutations to an element
for a class of oriented matroids. Of course, we have $L \ge 1$ for all Las Vergnas and Mandel oriented matroids and
in the realizable case there is the classical result that $L = \rk$ where $\rk$ is the rank of the oriented matroid. 
This was proved by Shannon in \cite{Shannon1979SimplicialCI}.
For connected Euclidean oriented matroids we show here that $L \ge 3$ if $\rk \ge 3$.
First, we assume that for each element $f$ exists a separable element $g$,
otherwise the oriented matroid is binary (hence realizable) or disconnected. We get $2$ mutations
on each 'side' of $f$ in the program $(\mathcal{O},g,f)$. We use then 
 the lexicographic extension to cut the two mutations in half.
Because the extended oriented matroid remains Euclidean, there must be a third one.

Generally, for Euclidean oriented matroids of rank $\rk$ it would be wishful 
to show $L = \rk$ because no counterexample is known. 
The case of rank $3$ was already shown by Levi in \cite{Levi26} (see \cite{OrientedMatroids}, Theorem 6.5.2 (i)).
But to obtain a proof analogous 
to Shannon's for the realizable case the assumption of Euclideaness seems to be not enough. 
Shannon uses for his proof the fact that we can always construct a hyperplane in general position going through $\rk -1$ points.
This translates to the intersection property $IP_1$ in the oriented matroid case 
introduced by Bachem/Wanka in \cite{BachemWanka}. We give a proof similar to Shannon's using the $IP_1$ and some subsequent perturbations to construct an element
'going through' $\rk -1$ mutations. Then we use again a lexicographic extension contravariant to that element to destroy all these mutations
and to find another one, like in the proof for $L \ge 3$.
To get the proof working, the extension yielded by the $IP_1$ must be Euclidean.

The intersection properties in \cite{BachemWanka} are not well-examined for oriented matroids,
even their equivalence in rank $4$ (which holds in the matroid case) can not be transferred to the oriented matroid case without again
assuming more, namely that we can choose the desired extensions to stay Euclidean. 
We know that Euclideaness is preserved if we extend a Euclidean oriented matroid program with a parallel extension (see \cite{OrientedMatroids}, Chapter 10.5),
but that does not need to hold for other oriented matroid programs,
we give a counterexample related to the Vamos matroid. In this paper we provide a short survey for what is known about preservation of Euclideaness and we show some more facts. It is preserved by direct sums, duality, minors, lexicographic extension, union, by exchanging target function and infinity
and by substituting the target function (or infinity) with an inseparable element.

From each mutation in an oriented matroid we can derive another oriented matroid (called a {\em mutant}) via a mutation-flip (see \cite{OrientedMatroids}, chapter 7.3).
We obtain a {\em mutation graph} if we take oriented matroids as vertices and draw edges between mutants. 
It is shown in \cite{Knauer} (with extensive computerhelp) that for uniform oriented matroids with $ \leq 9$ elements the mutation graph is connected.
We show here that an oriented matroid program $(\mathcal{O},g,f)$ stays Euclidean after a mutation-flip if either $f$ or $g$ (but not both) are adjacent to the mutation. 
All $8$-point non-Euclidean oriented matroids are counterexamples for the other cases, see also \cite{ExtensionSpaces} for more counterexamples.

From that, we can deduce that if a {\em totally non-Euclidean} oriented matroid (which means it has no Euclidean oriented matroid programs)
is connected to a Euclidean oriented matroid in the mutation-graph, the path in the graph must contain at least three edges (flips).
The $R(20)$ and its dual are the only known totally non-Euclidean oriented matroids (we show that in a paper appearing soon).
It would be helpful to find more examples or to check if a mutation-free element forces an oriented matroid to be totally non-Euclidean.
The inverse direction is false, the dual of $R(20)$ has no such element but is totally non-Euclidean. 
That shows also that $ \frak{O}_{Las-Vergnas}$ is not closed under duality, we do not know 
if this is true for $ \frak{O}_{Mandel}$ or not.

Finally, we apply this theory to show that uniform rank-$4$ oriented matroids are Mandel if they have a Euclidean mutant.
This holds also for uniform minimal non-Euclidean oriented matroids in all ranks. We take such an oriented matroid,
extend it lexicographically by cutting off a vertex from the mutation (simplicial tope) which was used by the flip. Now the mutation is adjacent to our new element and we flip
it back again. The new element gives then the extension desired for Mandelness.
It is mentioned in \cite{bokowski1990classification} (and easy to check with computerhelp)
that all uniform rank-4 oriented matroids with $8$ elements have a linear mutant.
Therefore, they are all Mandel which shows that Inclusion (3) is proper.
That means also the $RS(8)$, found independently by Roudneff and Sturmfels, see \cite{RoudneffSturmfels}, Chapter 5, is Mandel. It has rank $4$ but an element with only $3$ adjacent mutations, hence for the class of Mandel oriented matroids of rank $4$ holds $L < 4$.

The paper is build as follows. In the preliminary Section \ref{section:Preliminaries} we give the references to the literature needed to read this paper
and necessary definitions/special notations, if they differ from literature. 
As preliminary facts we show first that an oriented matroid in which each element has at least one corresponding
separable element must be non-binary or disconnected. But because all connected binary oriented matroids are realizable, 
we can omit them in our further investigations. Second, we recall some characterizations for simplicial topes/mutations
from literature, we did not find proofs for all of them, we close the gaps for completeness. Later, we will use mostly the characterization
of simplicial topes, that there is a base whose base-cocircuits are conformal. At the end of this section,
we show that disconnected realizable oriented matroids without coloops with $n$ elements have at least $3n - 9$ mutations.

In Section \ref{section:MandelnessAndMutations} we show that each element in an oriented matroid that is Mandel has always one adjacent mutation.
Hence each Mandel-oriented matroid is Las Vergnas. The proof is not very different to Mandel's proof in its dissertation but is using more modern language.
We characterize the property of Mandelness a bit more, the class of Mandel oriented matroids is nearly minor-closed.
Only the removal of so many elements that the rank is reduced may not preserve Mandelness.  
It is an open question if that class is closed under duality. 

In Section \ref{section:LexExtAndMutations} we examine the lexicographic extension further.
We show that it creates, destroys or preserves mutations like we mentioned before. 
The lexicographical extension
is in a way symmetrical. If we start with an oriented matroid $\mathcal{O} \cup f$, where $f$ is in general position, extend that lexicographically to 
$\mathcal{O} \cup \{f,f'\}$ with $f'$ being inseparable to $f$, remove $f$ to obtain $\mathcal{O} \cup f'$,
there is a (kind of reversed) lexicographic extension to obtain $\mathcal{O} \cup \{f,f'\}$ again.
This section is quite technically.

The next Section \ref{section:LasVergnasOrientedMatroids} gives some facts about Las Vergnas oriented matroids.
The class of those oriented matroids is not closed under duality and also not under taking minors,
we give a counterexample (the lexicographically extended $R(20)$). Moreover, it holds $L(\mathcal{O}_{LasVergnas}) = 1$.

The first main results come in Section \ref{section:EuclideanessAndMutations}. In Euclidean oriented matroids each element has three adjacent mutations,
we described the proof before. We use extensively the techniques of Section \ref{section:LexExtAndMutations}.
Then we give the analogon to Shannon's proof, but only for uniform rank $4$ oriented matroids having the (extended) $IP_1$-property.
They have at least $4$ adjacent mutations for each element.
Even though the proof follows Shannon's idea described before, it is rather techniqually and there is a lot of case checking to do.
Finally, we discuss the intersection properties further here.

Section \ref{section:Non-EuclideanOrientedMatroids} is a little excursus, we show that a directed cycle can not only contain the cocircuits of one simplicial tope
and we give a characterization for very strong components in minimal non-Euclidean oriented matorids. For every element $e$ of the groundset
there must be an edge containing $e$ and cocircuits having $e=+$ and $e=-$ in the component.

Section \ref{section:PreservationOfEuclideaness} is a short overview about preservation of Euclideaness.
It is known that it is preserved by duality, taking minors and Euclidean oriented matroid programs stay Euclidean by the extensions
they provide.  Other programs do not need to stay Euclidean, we give a counterexample.
We proved already in \cite{LexExt} that Euclideaness is preserved by lexicographic extension and by exchanging
target function and infinity. Here we show that it is preserved by direct sum and union, replacing target function (or infinity) by an inseparable element,
by a kind of perturbation process
and also by mutation flips if the target function is adjacent to the mutation but not the infinity. 
As an application we show, that in the rank $4$ case, if a totally Non-Euclidean oriented matroid is
connected to a Euclidean matroid in the mutation-graph, it is connected via at least three mutation flips. 

In Section \ref{section:NonEuclideanButMandel}
we show that a uniform rank $4$ or uniform minimal non-Euclidean rank $\rk$ oriented matroid is Mandel if it is connected via a mutation-flip to a Euclidean oriented matroid. For that we use again the fact that Euclideaness is preserved by these special mutation-flips and that we can exchange the order
of taking a mutation-flip and a special lexicographic extension. 

We finish in Section  \ref{section:concludingRemarks} with some concluding remarks.
We summarize what is now known about the minimal number of mutations adjacent to an element 
in the classes of oriented matroids we were considering here. Last not least, 
we give a bunch of questions raised by this paper and others that remained open
and could be maybe interesting for further research.

\section{Preliminaries}\label{section:Preliminaries}

We use the definitions and notions of \cite{OrientedMatroids}, especially from Chapter 10, 3 and 4.
We use many constructions of \cite{OrientedMatroids}, Chapter 7, and notions/results of  \cite{LexExt}, we will refer to it when needed.
If $\mathcal{O}$ is an oriented matroid we call its groundset $\mathcal{O}(E)$ (or sometimes only $E$) and if $g \in E$ we call $(\mathcal{O},g)$ 
the subset of the covectors of $\mathcal{O}$ having $g = +$, $S^0_g$ the pseudosphere of all covectors of  $\mathcal{O}$ having $g = 0$,
$S^{+0}_g$ the pseudosphere of all covectors of  $\mathcal{O}$ having $g \in \{0,+\}$ and 
$S^+_g$ (or $S^-_g$)  the pseudosphere of all covectors of  $\mathcal{O}$ having $g = +$ (or $g = -$), hence $S^+_g = (\mathcal{O},g)$.
If cocircuit elimination of an element $g$ between two cocircuits $X$ and $Y$ yields a unique cocircuit $Z$, we write sometimes {\em $El(X,Y,g) = Z$}
and we write $\Dir(X,Y)_{f,g} = Z_f$ like in \cite{LexExt}, Definition 3.1 for the direction of the edge $(X,Y)$ in the cocircuit graph of the oriented matroid program 
$(\mathcal{O},g,f)$.
We say a cocircuit $X$ is a {\em neighbour} of a cocircuit $Y$ if $X$ and $Y$ are conformal and comodular.
We start with some preliminary propositions. 

\begin{prop}\label{prop:GfNotEmpty}
Let $\mathcal{O}$ be an oriented matroid of rank $\rk$ with $n > \rk$ elements and $f \neq g  \in E$ both not being coloops and with $g$ lying in general position
and $f$ not being a loop.
Then there are cocircuits in $(\mathcal{O},g)$ with $f \neq 0$.
\end{prop}

\begin{proof}
Because $g$ is not a coloop, the oriented matroid $\mathcal{O}' = \mathcal{O} \setminus g$ remains of rank $\rk$.
Because $f$ is also not a loop in $\mathcal{O}'$, it can not lie in all hyperplanes of $\mathcal{M}(\mathcal{O}')$
(because all flats in a matroid are intersections of hyperplanes, see ...). 
Hence there is a cocircuit $X$ in $\mathcal{O}'$ with $f \notin z(X)$ hence $X_f \neq 0$
and because $g$ is in general position we have $g \notin z(X)$ hence $X_g \neq 0$ in $\mathcal{O}$.
\end{proof}

We cite three facts about of binary and regular matroids, see \cite{Oxley}, Theorem 6.5.4,  \cite{OrientedMatroids} Proposition 7.9.3 and \cite{OrientedMatroids}, Corollary 7.9.4. A matroid is called {\em binary} if it is representable over the field $GF(2)$ and {\em regular} if it is representable over every field. 

\begin{theorem}\label{theo:BinaryU24Minor}
A matroid is binary if and only if has no $U_{2,4}$ minor. 
\end{theorem}

We cite a very useful property for connected non-binary oriented matroids, see \cite{Oxley}, Proposition 12.3.7.
An oriented matroid is {\em connected} if its underlying matroid is {\em connected} (which is synonymous to {\em 2-connected}, see \cite{Oxley}, Proposition 4.1.3 and
the paragraph before). 

\begin{prop}\label{prop:U24MinorUsesElement}
Let $M$ be a connected matroid having a $U_{2,4}$-minor and suppose $e \in E(M)$. Then $M$ has a $U_{2,4}$-minor
containing $e$.
\end{prop}

\begin{theorem}\label{theo:BinaryOrientableIsRegular}
A binary matroid is orientable if and only if it is regular. 
\end{theorem}

\begin{theorem}\label{theo:OrientationsOfRegularMatroidsAreRealizable}
If $\mathcal{M}$ is a regular matroid, then all orientations of $\mathcal{M}$ are realizable. 
\end{theorem}

To the notion of inseparable elements, see \cite{OrientedMatroids}, Chapter 7.8.

\begin{prop}\label{prop:ExistNonInseparablePairs}
Let $\mathcal{O}$ be a simple connected non-binary oriented matroid and $f \in E$.
Then there is an element $g \in E \setminus f$ not inseparable to $f$.
\end{prop}

\begin{proof}
Let $(f,g)$ be a separable pair in a minor $\mathcal{O}'$ obtained from $\mathcal{O}$.
Then there are two cocircuits $X,Y$ in $\mathcal{O}'$ having $X_f = X_g$ and $Y_f = -Y_g$.
These cocircuits exist also in the oriented matroid $\mathcal{O}$ hence $f$ and $g$ are also separable there.
Hence we are done if we find a separable element to $f$ in a minor of $\mathcal{O}$.

Because $\mathcal{O}$ is non-binary, Theorem \ref{theo:BinaryU24Minor} yields that the corresponding matroid has a $U_{2,4}$-minor, 
and Proposition \ref{prop:U24MinorUsesElement} yields that it has a $U_{2,4}$-minor containing $f$,
hence a line $L$ with $4$-elements in $\mathcal{O}$. This is a rank $2$ oriented matroid with $4$ elements in general position. 
The other $3$ elements can not be all inseparable to $f$ and we are done.
\end{proof}
 
 \begin{lemma}\label{lem:ElementWithOnlyInseparableElIsRealizable}
 Let $\mathcal{O}$ be a simple connected oriented matroid.
 If there exist an element $f$ with all other elements $g$ being inseparable to $f$, then $\mathcal{O}$ is realizable.
 \end{lemma}

 \begin{proof}
 If there is such an element, then from Proposition \ref{prop:ExistNonInseparablePairs} before follows that $\mathcal{O}$ must be binary.
 From Theorem \ref{theo:BinaryOrientableIsRegular} follows that $\mathcal{O}$ must be regular
 and from Theorem \ref{theo:OrientationsOfRegularMatroidsAreRealizable}
 that $\mathcal{O}$ must be realizable then.
 \end{proof}
 
 We remark here that this lemma is not trivial at all. It uses the three theorems before which are far-reaching results of matroid theory.
 Now, we collect some well-known
 facts about mutations in oriented matroids.  Proofs for these results in literature were partly left as exercise for the reader,
we show here the missing gaps for completeness. 
For the notion of simplicial topes, mutations and mutation flips, we refer to \cite{OrientedMatroids}, Chapter 7.3.
We recall the definition of a {\em simplicial tope / mutation}.

\begin{definition}
A tope $T$ in an oriented matroid is called {\em simplicial} (or a {\em mutation}) if the lattice $[0,T]$ is boolean. 
\end{definition}

First, we need a lattice-theoretical result which is also known in literature but we did not find an explicit proof for it.
For the lattice theoretic notation we refer to \cite{OrientedMatroids}, Chapter 4.
Recall that a lattice is {\em relatively complemented} if in each interval $[x,y]$ and for an element $x < z < y$ there is another element $z'$ with 
$x < z' < y$ and $z \vee z' = y$ and $z \wedge z' = x$. We call it {\em atomic} (or {\em coatomic}) if each element is the join of atoms (or the meet of coatoms) and
we say an interval $I = [x,y]$ is $\vee$-regular and $\wedge$-regular if $y$ (resp. $x$) is the join (the meet) of all atoms (coatoms) of $I$.
The big face lattice $\mathcal{F}$ of an oriented matroid $\mathcal{O}$, hence also each tope lattice $\mathcal{T}$ of a tope of $\mathcal{O}$, has the 'diamond-property' which means that all intervals of length $2$ in $\mathcal{T}$ have cardinality $4$, this is \cite{OrientedMatroids}, Theorem 4.1.14 (ii).
Theorem 2 in \cite{BJORNER1981325} yields then that $\mathcal{T}$ is relatively complemented, atomic and coatomic
and that all intervals in $\mathcal{T}$ are $\vee$-regular and $\wedge$-regular.

\begin{prop}\label{prop:intervalsInTopeLattices}
An interval $I = [x,y]$ of length $n$ in a tope-lattice $\mathcal{T}$ of an oriented matroid $\mathcal{O}$ has at least $n$ atoms 
and $n$ coatoms.
\end{prop}

\begin{proof}
We can assume that $[x,y] = [0,T]$ is the tope lattice $\mathcal{T}$ itself, because of \cite{OrientedMatroids}, Proposition 4.1.9.
We show that the interval $[0,T]$ of length $n$ must have at least $n$ atoms.
We use induction, the length $2$ case is the diamond property. We assume the statement holds for $n-1$.
Let $x_{n-1}$ be a coatom in $[0,T]$. The induction assumption yields that the interval $[0, x_{n-1}]$ has at least $n-1$ atoms. Let $x_{n-2}$
be a coatom in $[0, x_{n-1}]$.
Because of the diamond property there is an element $z$ with $x_{n-1} \vee z = T$ and $x_{n-1} \wedge z = x_{n-2}$.
The interval $[0,z]$ has length $n-1$ hence has again at least $n-1$ atoms. One of these atoms cannot belong to the atoms in the interval
$[0,x_{n-1}]$ otherwise we would have $x_{n-1} = z$ because the intervals $[0,x_{n-1}]$ and $[0,z]$ are $\vee$-regular.
But all these atoms are covered by $T$ hence we have at least $n$ atoms in $[0,T]$. The case of coatoms works analogously.
\end{proof}

The if-direction of the next Proposition is obvious, the other direction needs more effort.

\begin{prop}\label{prop:mutationRCocircuits}
A tope $\mathcal{T}$ in an oriented matroid $\mathcal{O}$ of rank $\rk$ is simplicial iff it has only $\rk$ adjacent cocircuits. 
\end{prop}

%

\begin{proof}
For the if-direction recall
that a finite Boolean lattice is isomorphic to the lattice of a finite power set ordered by inclusion. Hence if the lattice has rank $\rk$ then it has $\rk$ atoms.
We show the other direction in four steps.

First, we show that in the tope lattice $\mathcal{T} = [0,T]$ we have exactly $\rk$ coatoms iff we have exactly $\rk$ atoms.
Proposition \ref{prop:intervalsInTopeLattices} yields that each coatom is the join of at least $\rk -1$ atoms. Hence if we have exactly $\rk$ atoms there are only $\rk$ combinations of $\rk -1$-sets of these atoms.
Hence we can only have $\rk$ cotaoms. The proof of the inverse direction is symmetric.

Second, we show that for each element $y$ of rank $n \le \rk$ the interval $[0,y]$ has $n$ atoms and coatoms.
We use induction going from $n$ to $n-1$. The base case $n=r$ and $y = T$ holds. 
We consider an interval $[0,x]$ of length $n-1$.
Let $y$ cover $x$. Because  $[0,y]$ has length $n$, the induction assumption yields that $[0,y]$ has $n$ atoms. 
Proposition \ref{prop:intervalsInTopeLattices} yields that the interval $[0,x]$ has at least $n-1$ atoms. 
If it would have $n$ atoms, these atoms would be also the atoms of $[0,y]$ and their join would be $x$, not $y$,
contradicting $[0,y]$ being $\vee$-regular.
Hence $[0,x]$ has exactly $n-1$ atoms (and coatoms which we proved before).
Then means also that each element $y$ of rank $n$ is the join of $n$ atoms because the interval $[0,y]$ is $\vee$-regular.


Third, we show that the join of $n$-atoms gives a unique element of rank $n$ for $n \le \rk$. The uniqueness follows from the associativity of the join. 
We assume there is a set $S$ of $n$-atoms whose join is an element $z$ of rank $r_z \neq n$.
 The interval $I = [0,z]$ has $r_z$ atoms. 
  Each atom of $S$ lies in $I$ hence  $r_z < n$ cannot appear.
If $r_z > n$ we would obtain a set of $r_z - 1$ atoms of $I$  whose join is $z$ and not a cotaom of $I$.
But $I$ has  $r_z$ coatoms which are the join of $r_z -1$ atoms, hence we would not have enough sets of $r_z - 1$ atoms of $I$ to obtain all these coatoms.
We obtain $r_z = n$.
We proved here that the lattice is modular. We have 
\[ \rk(x) + \rk(y) - \rk(x \wedge y)= \rk(x \vee y). \]
Fourth, we can identify the elements $x$ of $\mathcal{T}$ with the subsets $S(x) \in \mathcal{P}(\{p_1, \hdots, p_{\rk}\})$ of atoms whose join they are.
It is clear that we have $S(x \wedge y) = S(x) \cap S(y)$ and from modularity follows also $S(x \vee y) = S(x) \cup S(y)$.
The lattice $\mathcal{T}$ is Boolean.
\end{proof}

The next Lemma is an important characterization of simplicial topes, see \cite{OrientedMatroids}, Proposition 7.3.7. 

\begin{lemma}\label{prop:BasicCocircuitsOfASimplicialTope}
 Let $\mathcal{T}$ a positive tope in an oriented matroid of rank $\rk$. Then $\mathcal{T}$ is simplicial iff there is a base 
 $B = \{ b_1, \hdots, b_{\rk} \}$ such that the basic cocircuits $c^*(b_i, B)$ are all positive. In that case
 in the corresponding pseudosphere arrangement the tope $\mathcal{T}$ is bounded by the pseudospheres $S^{+0}_{b_i}, \hdots, S^{+0}_{b_{\rk}}$. 
 Or each $b_i$ is lying in the complement of exactly one subtope of $\mathcal{T}$. 
 We say the $b_i$ are {\em adjacent} to the mutation and write sometimes for the mutation $M = [b_1, \hdots, b_{\rk}]$ and we call the basic cocircuit $c^*(b_i, B)$
 sometimes the {\em $b_i$-base cocircuit of $M$}.
\end{lemma}

\begin{proof}
In \cite{MandelDiss}, Chapter VI, Proposition 5, page 320, this is the equivalence from (5.1) to (5.4). The proof of the equivalences was left as an exercise there.
The equivalences from (5.2), (5.3) and (5.4) were proved by Las Vergnas
in \cite{LasVergnas}, Section 1. The equivalence from (5.1) to (5.3) is our Proposition \ref{prop:mutationRCocircuits}. That gives the missing links.
\end{proof}

\begin{remark}\label{rem:BasicCocircuitsOfASimplicialTope}
If $\mathcal{T}$ is not positive we may assume all base cocircuits being pairwise conformal (instead of being positive) and their composition yields $\mathcal{T}$.
\end{remark}

Finally, we show some facts about the number of mutations in disconnected oriented matroids. 
Recall that in a connected realizable oriented matroid $\mathcal{O}$ of rank $\rk$ and $n$ elements each element has at least $\rk$ adjacent mutations
and  $\mathcal{O}$ has at least $n$ mutations, see the introduction and Section \ref{section:EuclideanessAndMutations} for further investigations. See \cite{OrientedMatroids}, Chapter 7.6, for the definition of the direct sum of two oriented matroids.

\begin{lemma}
Let $\mathcal{O}$ be the direct sum of two matroids $\mathcal{O}^1$ and $\mathcal{O}^2$.
Let $f \in \mathcal{O}^1(E)$ and $g \in \mathcal{O}^2(E)$.
If $f$ has $m$ adjacent mutations  in $\mathcal{O}^1$ and if $\mathcal{O}^2$ has $n$ mutations
then $f$ has $m*n$ adjacent mutations in $\mathcal{O}$.
\end{lemma}

\begin{proof}
The direct sum of two mutations stays a mutation. 
\end{proof}

\begin{lemma}\label{lem:disconnectedOMNumberMutations}
A simple disconnected realizable oriented matroid $\mathcal{O}$ without coloops of rank $\rk$ with $n$ elements has at least  
$3n - 9$ mutations.
\end{lemma}

\begin{proof}
Because $\mathcal{O}$ has no coloops, it is the direct sum of at least $2$ components of rank $\ge 2$ hence it must have at least rank $4$ and
$n \ge 6$ elements. Let  $\mathcal{O} = \mathcal{O}^1 \times \mathcal{O}^2$ be the direct sum
of two connected components $\mathcal{O}^1$ and $\mathcal{O}^2$ of rank $\rk_1$ and $\rk_2$ with $n_1$ and $n_2$ elements, they have
at least $n_1$ and $n_2$ mutations because they are realizable.
We have $\rk = \rk_1 + \rk_2$ and $n = n_1 + n_2$.
For $i \neq j \in \{1,2\}$ each element in $O^i(E)$ has at least $\rk_i * n_j$ adjacent mutations. 
We count the adjacent mutations for each element and divide through $\rk$ because each mutation was counted $\rk$ times.
We obtain $(n_1 * \rk_1 * n_2 + n_2 * \rk_2 * n_1) / \rk = n_1 * n_2$ mutations for $\mathcal{O}$.
It must hold $n_1,n_2 \ge 3$ because we have no coloops.
We have to choose $n_1,n_2$ such that $n_1 * n_2$ is minimal.
Let $n_1 =  \lfloor n/2 \rfloor - x$ and $n_2 =  \lfloor (n+1)/2 \rfloor + x$. Then $n_1 * n_2$ is minimal iff $x$ is maximal,
hence we choose $n_1 = 3$ and $n_2 = n - 3$. 
We obtain $3n - 9$ mutations which is more than $n$ because $n \ge 6$. 
The estimation gets worse if we have more than $2$ connected components.
\end{proof}

\section{Mandel Oriented Matroids and mutations}\label{section:MandelnessAndMutations}

The fact, that Euclideaness implies the existence of simplicial topes in an oriented matroid was already mentioned in Mandels 
dissertation \cite{MandelDiss} and in \cite{Knauer}, Theorem 4.5. We reformulate this result in another way.
The next Theorem is similar to Theorem (VI), III in \cite{MandelDiss}, page 318.

\begin{theorem}\label{theo:simplicialTopeMandel}
Let $(\mathcal{O},g,f)$ be a Euclidean oriented matroid program of rank $\rk$ with $g$ lying in general position and not being a coloop
and $f$ being not a loop or a coloop.
Let $(\mathcal{O} / f,g,e)$ be Euclidean oriented matroid programs for all $e \notin \{f,g\}$ that are not loops/coloops.
Let $G_f$ be the cociruit graph of the program $(\mathcal{O},g,f)$.
Then the graph $G^+_f$ (or $G^-_f$) of all cocircuits with $f = +$ (or $f = -$)  of $G_f$ is not empty
and there is a simplicial tope $T$ in $(\mathcal{O},g,f)$ adjacent to $S_f$ in the pseudosphere arrangement and with all cocircuits of $[0,T]$ having 
$g = +$ and $f = 0$ 
except for one cocircuit $X$ having $f = +$ (or $f = -$) and $g = +$.
That cocircuit $X$ has no predecessors (or successors) in the cocircuit graph $G^+_f$ (or $G^-_f$).
Each cocircuit without predecesors yields such a simplicial tope.
\end{theorem}

\begin{proof}
We prove the statement by induction over the rank of the oriented matroid. The rank 2 case is obvious. We take a cocircuit
$X$ with $X_f = 0$, it has $X_g \neq 0$ because $g$ is in general position. Because we have at least $3$ elements one of the two neighbours 
$Y$ of $X$ has $Y_g \neq 0$ and $T = X \circ Y$ is the desired simplicial tope.
We assume the Theorem holds for all oriented matroids of rank $\rk-1$.
Because $(\mathcal{O},g,f)$ is Euclidean the graph $G_f$ has no directed cycles. Proposition \ref{prop:GfNotEmpty} yields that
 $G_f$ (we may assume $G_f^+$) is not empty. Then there is at least one source in $G^+_f$ which is
a cocircuit $X$ with $X_f = X_g = +$ without a predecessor in $G^+_f$.
Using reorientation which does not affect our argumentation we may assume that $X$ is a positive cocircuit.

Now let $X$ be lying on an edge $F$. Then, because $X$ has no predecessor in $G^+_f$, there must lie a neighbour $Y$ of $X$ in $F$
not being in $G^+_f$, hence having $f = 0$ or $g = 0$. 
The case that $Y$ has $Y_f = Y_g = 0$ can not appear
because $g$ is in general position.
In the case $Y_g = 0$ and $Y_f = +$ cocircuit elimination of $f$ between $-Y$ and $X$ yields
$Z$ with $Z_f = 0$ and $Z_g = +$, the edge is directed from $Z$ to $Y$ via $X$. But because $X$
has no predecessor in $G^+_f$ the cocircuit $Z$ must also be a neighbour of $X$.
Hence we obtain always a neighbour of $X$ lying in $F$ having $f = 0$ and $g = +$.
Because of $\rk \ge 3$ and Proposition \ref{prop:intervalsInTopeLattices}, there are two different edges going through the cocircuit $X$.
Let $Y^1$ and $Y^2$ be the two neighbours of $X$ lying on these two edges having $f = 0$ and $g = +$.
It holds $Y^1 \neq Y^2$.

The oriented matroid program $\mathcal{O}' = (\mathcal{O} \setminus \supp(X)) \cup \{f,g\} / f$ has rank $\rk - 1$.
Then $Y'^{1(2)} = Y^{1(2)} \setminus \supp(X)$ are different cocircuits in $(\mathcal{O}',g)$ hence $g$ is not a coloop and is in general position in $\mathcal{O}'$.
On the other hand let $Y'$ be a cocircuit of $(\mathcal{O}',g)$ hence $\rk(z(Y')) = \rk-2$ in $\mathcal{O}'$. 
Let $Y$ be the extended cocircuit of $Y'$ in $\mathcal{O}$ (hence $Y_f = 0$). It holds 
$z(X \circ Y) \supseteq z(Y')$ hence $\rk(X \circ Y) = \rk-2$ in $\mathcal{O}$.
We have $\sep(Y,X) = \emptyset$ otherwise cocircuit elimination would yield a cocircuit $Z$ 
between $Y$ and $X$ with $Z_f = +$. We would have $\Dir(Z,X) = \Dir(Y,X) = +$ hence $Z$ would be
a predecessor of $X$ which is impossible.

Let $e \in z(X \circ Y^1) \setminus z(X \circ Y^2)$. Then
$(\mathcal{O}',g,e)$ is a Euclidean oriented matroid program because of the asssumptions
with $e \neq g$ not being a loop and $g$ lies in general position and is not a coloop in $\mathcal{O}'$.
We have $Y'^1_e = 0, Y'^2_e \neq 0$ and $Y'^1_g = Y'^2_g = +$ in $\mathcal{O}'$.
Then cocircuit elimination of $g$ between $-Y'^1$ and $Y'^2$ yields a cocircuit $Z'$ with $Z'_g = 0, Z'_e = Y'^2_e \neq 0$. We obtain two 
cocircuits in $\mathcal{O}'$ having $e \neq 0$, hence $e$ is not a coloop. 
The oriented matroid $\mathcal{O}'$ has $\ge \rk$ elements and is of rank $\rk-1$ hence Proposition \ref{prop:GfNotEmpty} holds
and either $G^+_e$ or $G^-_e$ is not empty. 

Because of the induction assumption it has a simplicial tope $T'$ 
having $g= +$ for all cocircuits of $[0,T']$.
Again using reorientation does not affect anything of that proof: because $X_e = 0$ for all $e \in \mathcal{O}' \setminus \{f,g\}$ we may assume $T'$ being positive. Because of Lemma \ref{prop:BasicCocircuitsOfASimplicialTope}, $\mathcal{O}'$ has a base $B' = \{b_1, \hdots, b_{\rk - 1}\}$
with $g \notin B'$ such that the cocircuits $Y'^i = c^*(b_i,B')$ are all positive.
The extended cocircuits $Y^i$ are neighbours to $X$ and have $Y^i_f = 0$, they stay positive. 
Then $B = \{ b_1, \hdots, b_{\rk - 1}, f\}$ is a base with $Y^i = c^*(b_i,B)$ because $Y^i_f = 0$.
On the other hand we have $X_{b_i} = 0$ because $g \notin B'$ hence we have $X = c^*(f,B)$.  
We get a simplicial tope $T$ in $\mathcal{O}$ adjacent to $S_f$ and with all cocircuits of $[0,T]$ having $g = +$ because all cocircuits 
$Y^i$ have $g = +$ and $X_g = +$ and $T = Y^1 \circ \hdots \circ Y^{\rk - 1} \circ X$.
%
\end{proof}

The next theorem is also proved in \cite{Knauer}, Corollary 4.6.
We call an oriented matroid {\em Las Vergnas} if each element of its groundset has an adjacent mutation.

\begin{theorem}\label{theo:MandelIsLAsVergans}
Let $\mathcal{O}$ be a Mandel oriented matroid. Then each element $f$ not being a loop/coloop in $\mathcal{O}$ has an adjacent mutation, hence $\mathcal{O}$ is a Las-Vergnas oriented matroid.
\end{theorem}

\begin{proof}
Let $\mathcal{O}' = \mathcal{O}\cup g$ the extension in general position (and not being a coloop) such that all programs $(\mathcal{O},f,g)$ are Euclidean for all elements $f$ of the groundset of $\mathcal{O}$ that are not loops/coloops. Then also the programs  $(\mathcal{O},g, f)$ are Euclidean and (because Euclideaness is minor-closed) the programs $(\mathcal{O} / f,g,e)$ for all $e \notin \{f,g\}$ not being loops/coloops. These are exactly the assumptions of Theorem \ref{theo:simplicialTopeMandel}.
\end{proof}

As a corollary we obtain three examples of non-Mandel oriented matroids because they are not Las-Vergnas.

\begin{corollary}
The oriented matroids from \cite{BokowskiMutationProblem}, Chapter 6, \cite{hall2004counterexamples}, Chapter 2.4 and $R(20)$, given in \cite{Richter-Gebert}, Chapter 2, are not Mandel.
\end{corollary}

We continue characterizing Mandel oriented matroids.

\begin{theorem}\label{theo:MandelOMsAreMinorClosed}
A Mandel oriented matroid remains Mandel after a contraction and/or after deletion of elements if that deletion does not reduce the rank.
\end{theorem}

\begin{proof}
Let $\mathcal{O}$ be a Mandel oriented matroid and let $\mathcal{O}' = \mathcal{O} \cup f$ the extension in general position (not being a coloop)
such that all programs $(\mathcal{O}', g, f)$ are Euclidean for $g$ not being loops/coloops. Then in the programs $(\mathcal{O}' \setminus e, g, f)$ and 
$(\mathcal{O}' / e, g, f)$ the element  $f$ is also in general position. These programs stay Euclidean if $f$ is not a coloop there.
The coloop-case can only appear in a deletion of elements that reduces the rank of $\mathcal{O}$.
\end{proof}

Because all rank $3$ oriented matroids are Euclidean we have: 

\begin{theorem}\label{theo:MandelOMsAreMinorClosed}
The class of Mandel oriented matroids with rank $\le 4$ is minor closed.
\end{theorem}

For uniform oriented matroids holds:

\begin{theorem}\label{theo:MandelOMsAreMinorClosed}
The class of uniform Mandel oriented matroids is minor closed.
\end{theorem}

\begin{proof}
We have to exclude the coloop-case here. If we remove elements of a uniform oriented matroid of rank $\rk$ until the rank is reduced it remains
an oriented matroid of rank $\rk -1$ with $\rk -1$ elements which is too small to be considered.
\end{proof}

We add a small Remark here.

\begin{remark}\label{rem:dualOfR20MandelOrNot}
We do not know if the class of Mandel oriented matroids is closed under duality. If it were closed it would be also minor-closed in general.
The dual of $R(20)$ has no mutation-free element (this is yielded by computerhelp), hence it is Las Vergnas, but we do not know if it is Mandel.
\end{remark}

%

\section{the lexicographic extension and its mutations}\label{section:LexExtAndMutations}

For the definition of lexicographic extension, see \cite{OrientedMatroids}, Chapter 7.2, Definition 7.2.3 and Proposition 7.2.4.
We show in this Section, that the lexicographic extension destroys mutations adjacent to the element $f$ where it is inseparable to, creates at least one new adjacent mutation and leaves non-adjacent mutations as they are. The new element $f'$ is in a way exchangeable to $f$ if $f$ lies in general position.
We need some preliminary Propositions.

\begin{prop}\label{prop:ExtensionInGeneralPositionSameNumberOfCocircuits}
Let $\mathcal{O} \cup f$ and $\mathcal{O} \cup f'$ two single-element extensions of an oriented matroid $\mathcal{O}$
with $f$ and $f'$ being in general position. Then both extensions have the same number of cocircuits.
\end{prop}

\begin{proof}
Because both extensions are in general position, we have only new cocircuits derived from an edges of $\mathcal{O}$.
If $F$ is an edge in $\mathcal{O}$ then $z(F) \cup f$ and $z(F) \cup f'$ are hyperplanes because $f$ and $f'$ are in general position
yielding cocircuits. Hence each edge in $\mathcal{O}$ yields two cocircuits in the extensions and other new cocircuits cannot appear.
\end{proof}

\begin{prop}\label{prop:InseparablePairsAndNieghbours}
Let $\mathcal{O}$ be an oriented matroid and let $(f,f')$ be an inseparable pair of elements in $E$ both being in general position. 
If $X$ is a cocircuit in $\mathcal{O}'$ with $X_{f'} = 0$ and $X_f \neq 0$ (or vice versa),
then we have a neighboured cocircuit $Y$ to $X$ having $Y_{f} = 0$ and $X_{f'} \neq 0$ (or vice versa).
It holds $X_f = - \alpha Y_{f'}$ (or $X_{f'} = - \alpha Y_f$) where $\alpha = +$ ($\alpha = -$) if $(f,f')$ is a contravariant (covariant) pair.
We call $Y$ the {\em corresponding cocircuit to $X$}. 
\end{prop}

\begin{proof}
We assume $f,f'$ being contravariant.
Let $F = z(X) \setminus f'$. Because $f'$ lies in general position, $F$ is an edge. We have $f \notin z(F)$.
Let $Y,Z$ be the two neighbours of $X$ lying on $F$. We have $Y_{f'} = -Z_{f'} \neq 0$
and if they had both $f \neq 0$ we would obtain $Y_f = Z_f \neq 0$ contradicting $f,f'$ being inseparable.
 We assume $Y_f = 0$.
If we had $Y_{f'} = X_f$, then $Z_f = X_f$ and $Z_{f'} = - X_f$ would hold contradicting $f,f'$ being contravariant.
From $z(Y) \setminus f = F$ (because $f$ is in general position) and $\sep(X,Y) = \emptyset$ follows $Y_e = X_e$ for all $e \in E \setminus f$.
\end{proof}

\begin{prop}\label{prop:LexExtensionCocircuitsHavingFPrimeZeroOrFZeroHaveNeighbours}
Let $\mathcal{O}$ be an oriented matroid of rank $\rk$ with independent elements $I = \{f,e_2, \hdots, e_r\}$ of the groundset $E$ and $f$ being in general position.
Let $\mathcal{O}' = \mathcal{O}[f^{\alpha_1},e_2^{\alpha_2}, \hdots, e_r^{\alpha_r}]$ be a lexicographic extension of $\mathcal{O}$.
If $X$ is a cocircuit in $\mathcal{O}'$ with $X_{f'} = 0$ and $X_f \neq 0$ (or vice versa),
then we have a neighboured cocircuit $Y$ to $X$ having $Y_{f} = 0$ and $X_{f'} \neq 0$ (or vice versa) and $Y_e = X_e$ for all $e \in E \setminus f$.
We have $X_f = -\alpha_1 Y_{f'}$ (or $X_{f'} = -\alpha_1 Y_f$).
\end{prop}

\begin{proof}
The elements $f,f'$ lie both in general position and are contravariant (covariant) if $\alpha_1 = +$ (or $\alpha_1 = -$).
\end{proof}

\begin{prop}\label{prop:LexExtensionCocircuitsHavingFPrimeZero}
Let $\mathcal{O}$ be an oriented matroid of rank $\rk$ with independent elements $I = \{f,e_2, \hdots, e_r\}$ of the groundset $E$ and $f$ being in general position.
Let $\mathcal{O}' = \mathcal{O}[f^{\alpha_1},e_2^{\alpha_2}, \hdots, e_r^{\alpha_r}]$ be a lexicographic extension of $\mathcal{O}$.
Then if $X$ is a cocircuit in $\mathcal{O}'$ having $X_{f'} = 0$ and $X_f \neq 0$ we have
$X_f = - \alpha_1* \alpha_2 X_{e_i}$ where $i$ is the first index of $I \setminus f$ such that $X_{e_i} \neq 0$.
\end{prop}

\begin{proof}
Proposition \ref{prop:LexExtensionCocircuitsHavingFPrimeZeroOrFZeroHaveNeighbours} yields a cocircuit $Y$ corresponding to $X$ with $Y_f = 0$.
We have 
\[X_f \overset{\text{Prop \ref{prop:LexExtensionCocircuitsHavingFPrimeZeroOrFZeroHaveNeighbours}}}{=} -\alpha_1 Y_{f'} \overset{\text{lex. Extension}}{=} -\alpha_1 * \alpha_2 Y_{e_i} \overset{\text{Prop \ref{prop:LexExtensionCocircuitsHavingFPrimeZeroOrFZeroHaveNeighbours}}}{=} - \alpha_1* \alpha_2 X_{e_i}.\]  
\end{proof}

The next Lemma shows that the lexicographic extension is in a way symmetric if the element $f$ where it is inseparable to lies in general position.

\begin{lemma}\label{lemma:LexExtensionIsomorphism}
Let $\mathcal{O}$ be an oriented matroid of rank $\rk$ with independent elements $I = \{f,e_2, \hdots, e_r\}$ of the groundset $E$ and $f$ being in general position.
Let $\mathcal{O}^2 = \mathcal{O} \cup f' = \mathcal{O}[f^{\alpha_1},e_2^{\alpha_2}, \hdots, e_r^{\alpha_r}]$ be a lexicographic extension of $\mathcal{O}$.
Let $\mathcal{O}^3 = \mathcal{O} \cup f' = \mathcal{O}[f^{\alpha_1},e_2^{-\alpha_1 * \alpha_2}, \hdots, e_r^{-\alpha_1 * \alpha_r}]$ be a lexicographic extension of $\mathcal{O}$.
Then $\mathcal{O}^3$ is isomorphic to $\mathcal{O}^2$. The isomorphism maps $f$ to $f'$ and vice versa and is the identity 
for all other elements of $E \setminus f$.
\end{lemma}

\begin{proof}
Let $X$ be a cocircuit in $\mathcal{O}^2$ with $X_f \neq 0$ and $X_{f'} \neq 0$. Then it holds $X_{f'} = X_f$ in  $\mathcal{O}^2$ and  $\mathcal{O}^3$.
The cocircuit $X$ stays the same in $\mathcal{O}^2$ and $\mathcal{O}^3$.
Let $Y$ be a cocircuit in $\mathcal{O}^2$ with $Y_f = 0$ and $Y_{f'} \neq 0$. Then we have 
$Y_{f'} = \alpha_2 Y_{e_i} \text{ in } \mathcal{O}^2$.
This cocircuit exist also in $\mathcal{O}^3$. We have there
$Y_{f'} = - \alpha_1 \alpha_2 Y'_{e_2} \text{ in } \mathcal{O}^3$.
Let $Y'$ be the corresponding cocircuit to $Y$ in $\mathcal{O}^3$ from Proposition \ref{prop:LexExtensionCocircuitsHavingFPrimeZeroOrFZeroHaveNeighbours}. 
It has 
$Y'_{f'} = 0, Y'_f = - \alpha_1 Y_{f'} =\alpha_2 Y_{e_i} \text{ in } \mathcal{O}^3 \text{ and }Y'_e = Y_e \text{ for all }e \in E \setminus f.$
Hence, if we exchange the values of $f$ and $f'$, then $Y'$ has the same values in $\mathcal{O}^3$ as $Y$ in $\mathcal{O}^2$.

Analogously if we have a cocircuit $Z$ with $Z_{f'} = 0$ and $Z_f \neq 0$ in $\mathcal{O}^2$, Proposition \ref{prop:LexExtensionCocircuitsHavingFPrimeZero} yields $Z_f = - \alpha_1 \alpha_2 Z_{e_i}$. Proposition \ref{prop:LexExtensionCocircuitsHavingFPrimeZeroOrFZeroHaveNeighbours}
 yields a corresponding cocircuit $Z'$ in 
$\mathcal{O}^2$
with $Z'_f = 0$ and $Z'_e = Z_e$ for all $e \in E \setminus \{f,f'\}$.
This cocircuit exists also in $\mathcal{O}^3$ and it holds there $Z'_{f'} = - \alpha_1 \alpha_2 Z'_{e_i}$.
Hence if we exchange $f$ and $f'$ the cocircuit $Z'$ has in $\mathcal{O}^3$ the same values as in $\mathcal{O}^2$.

Let $X$ be a cocircuit in $\mathcal{O}^2$ with $X_f = X_{f'} = 0$.
Then there are two neighbours $Y,Z$ in $\mathcal{O}$ with $\sep(Y,Z) = f'$ in $\mathcal{O}^2$ because of \cite{OrientedMatroids}, Proposition 7.1.4,
and with $Y_f = Z_f = 0$ because $f'$ is in general position.
From the paragraph before we obtain cocircuits $Y',Z'$ in $\mathcal{O}^3$ having $Y'_{f'} = Z'_{f'} = 0$, $\sep(Y',Z') = f$, $Y'_e = Y_e$ and $Z'_e = Z_e$ for all $e \in E \setminus f$.
Cocircuit elimination of $f$ between $Y',Z'$ yields a cocircuit $X'$ having $X'_e = X_e$ for all $e \in E \setminus \{f,f'\}$ and $X'_{f'} = X_f = 0$.

To each cocircuit $X$ in $\mathcal{O}^2$ we obtain a cocircuit $X'$ in $\mathcal{O}^3$ with the same values if we exchange $f$ and $f'$.
 This proves everything because $\mathcal{O}^2$ and $\mathcal{O}^3$ have the same number of cocircuits which is Proposition
\ref{prop:ExtensionInGeneralPositionSameNumberOfCocircuits}.
\end{proof}

We show that the lexicographic extension yields a mutation.

\begin{lemma}\label{lemma:LexExtensionYieldsOneMutation}
Let $\mathcal{O}$ be a uniform oriented matroid of rank $\rk$ with elements $[f,e_1,\hdots$ $\hdots,e_{r-1}]$ of the groundset.
Let $\mathcal{O}' = \mathcal{O}[f^+,e_1^+, \hdots,e_{r-1}^+]$ be a lexicographic extension of $\mathcal{O}$.
Then $M = [f,f',e_1,\hdots, e_{r-2}]$ is a mutation.
\end{lemma}

\begin{proof}
Let $X,Y,Z^1, \hdots, Z^{r-2}$ be the (pairwise comodular) base-cocircuits with 
\begin{equation*}
\begin{split}
z(X) &= \{f', e_1, \hdots, e_{r-2} \} \text{ with } X_f = +, \\ 
 z(Y) &= \{f, e_1, \hdots, e_{r-2} \} \text{ with } Y_{e_{r-1}} = -, \\
 z(Z^i) &= \{f, f', e_1, \hdots, e_{i-1}, e_{i+1}, \hdots e_{r-2} \} \text{ with } Z^i_{e_i} = +.
\end{split}
\end{equation*} 
The lexicographic extension yields $Y_{f'}= -$.
We show that all cocircuits are conformal.
Proposition \ref{prop:BasicCocircuitsOfASimplicialTope} and Remark \ref{rem:BasicCocircuitsOfASimplicialTope} yield a mutation then.
\begin{enumerate}
\item case $Y, Z^i$: We assume $h \in \sep(Y,Z^i)$. Then $El(Y,Z^i,h) = C$ with $C_h = C_f = C_{e_1}= \hdots = C_{e_{i-1}} = 0, C_{f'} = -$ and $C_{e_i} = +$ which is impossible because the lexicographic extension yields $C_{f'} = C_{e_i}$.
\item case $X,Z^i$: We assume $h \in \sep(X,Z^i)$. Then $El(X,Z^i,h) = C$ with $C_h = C_{f'} = C_{e_1}= \hdots = C_{e_{i-1}} = 0, C_{e_i} = C_{f} = +$
which is impossible because Proposition \ref{prop:LexExtensionCocircuitsHavingFPrimeZero} yields $C_f = - C_{e_i}$.
\item case $X,Y$: We assume $h \in \sep(Y,Z)$. Then $El(Y,Z,h) = C$ with $C_h = C_{e_1}= \hdots = C_{e_{r-2}} = 0, C_{f'} = -,C_f = +$
which is impossible because $f$ and $f'$ are contravariant. 
\item case $Z^i,Z^j$ with $i \neq j$: We assume $h \in \sep(Z^i,Z^j)$. 
This would yield $Z^i_h = -Z^j_h \neq 0$, but since $\mathcal{O}$ is uniform, we have $X_h \neq 0$, hence $Z^i_h = X_h$ and $X_h = Z^j_h$ which we proved before yielding a contradiction. \qedhere
\end{enumerate}
\end{proof}

We need a small Proposition for the next result.

\begin{prop}\label{prop:ElementInGeneralPositionInvolvedInMutation}
Let $\mathcal{O}$ be an oriented matroid of rank $\rk$ with an element $f$ of the groundset lying in general position.
Let $M$ be a mutation with a cocircuit $X$ having $X_f = 0$. Then $f$ is adjacent to the mutation.
\end{prop}

\begin{proof}
Let $M = (e_1, \hdots, e_r)$ be the mutation. $M$ has only adjacent base-cocircuits, hence
for $\rk -1$ elements $e_i$ of the mutation holds $X_{e_i} = 0$. 
Hence, if $f$ is not in $M$, in the oriented matroid $\mathcal{O} \setminus f$ the cocircuit $X$ remains a cocircuit
contradicting $f$ being in general position.
\end{proof}

The lexicographic extension leaves non-adjacent mutations as they are.

\begin{lemma}\label{lemma:LexExtensionleavesNonAdjacentMutations}
Let $\mathcal{O}$ be a uniform oriented matroid of rank $\rk$ with elements $\{f,e_2,\hdots ,e_{\rk}\}$ of the groundset.
Let $\mathcal{O}' = \mathcal{O}[f^{\alpha_1},e_2^{\alpha_2}, \hdots,e_{\rk}^{\alpha_{\rk}}]$ be a lexicographic extension of $\mathcal{O}$.
Let $M$ be a mutation in $\mathcal{O}$ not adjacent to $f$. Then $M$ is also a mutation in $\mathcal{O}'$. 
\end{lemma}

\begin{proof}
First, because we are in the uniform case, if we have a cocircuit of the mutation having $f = 0$ then $f$ must be adjacent to the mutation because of Proposition 
\ref{prop:ElementInGeneralPositionInvolvedInMutation}.
Hence all cocircuits of $M$ must have the same $f$-value $\neq 0$.
But then they have the same $f'$-value $\neq 0$ and no separating elements in $\mathcal{O}'$ are between them, we are done.
\end{proof}

Now we show that the lexicographic extension destroys mutations or 'shifts' its positions in a way. 

\begin{lemma}\label{lemma:MutationsAndLexExtensions}
Let $\mathcal{O}$ be a uniform oriented matroid of rank $\rk$ and let $M = (f,e_2, \hdots, e_r)$ be a mutation
in $\mathcal{O}$. Let the elements $\{f,e_2, \hdots, e_r\}$ be oriented versus the mutation.
Let $\mathcal{O}' = \mathcal{O}[f^+,g^-, \hdots]$ be a lexicographic extension of $\mathcal{O}$ with $g \notin \{f,e_2, \hdots, e_r\}$ and let the cocircuit $Y$ of the mutation with $Y_f = +$ have also $Y_g = +$. Then $M' = (f',e_2, \hdots, e_r)$ is a mutation in $\mathcal{O}'$ and the tope $T$ with $T_f = T_{e_2} =  \hdots =  T_{e_r} = +$ is not a simplicial region anymore. $f$ has no adjacent mutations with cocircuits having $f = +$ and $g = +$ where $f'$ is not involved. 
\end{lemma}

\begin{proof}
 First, we show that $M'$ is a mutation.
We look at the cocircuits $X^1, \hdots, X^{r-1}$ and $Y$ of $M$ having $X^1_f = \hdots = X^{r-1}_f = 0$ and $Y_f = +$ and
having $X^1_{e_2} =  \hdots = X^{r-1}_{e_r} = Y_f = +$ (it holds $X^i_{e_j} = 0$ for $j \neq i+1$ and all cocircuits have $g= +$). Then  
we have $Y_{f'} = +$, too. We have $X^1_{f'} = \hdots = X^{r-1}_{f'} = -$ and obtain new cocircuits $X'^1, \hdots, X'^{r-1}$ by elimination of $f'$ between $X^1,\hdots,X^{r-1}$ and $Y$.
The new cocircuits and $Y$ are pairwise neighbours and they are the base cocircuits of $M'$
 hence $M'$ is because of Proposition \ref{prop:BasicCocircuitsOfASimplicialTope} a mutation consisting of these cocircuits. The tope $T$ consists of the cocircuits
$X^1, \hdots, X^{r-1}, X'^1, \hdots, X'^{r-1}$, hence more than $\rk$ cocircuits. It is not a mutation anymore because of Proposition \ref{prop:mutationRCocircuits}. The rest follows immediately.
\end{proof}

We prove here more general some Propositions about mutations and contravariant pairs.

\begin{prop}\label{prop:MutationsForContravariatPairs}
Let $\mathcal{O}$ be an oriented matroid and let $(f,f')$ be a contravariant pair of elements.
Let $f'$ be in general position.
Let $M$ be a mutation adjacent to $f$ and let $X$ be the cocircuit of $M$ having $X_f \neq 0$.
Then if $f'$ is not involved in the mutation we have $X_{f'} = X_f$ and all cocircuits of the mutation have the same $f'$-value. 
\end{prop}

\begin{proof}
If there is a cocircuit of the mutation having $f' = 0$ then Proposition \ref{prop:ElementInGeneralPositionInvolvedInMutation}
 yields that $f'$ is involved in the mutation.
Hence if $f'$ is not involved all cocircuits of the mutation have $f' \neq 0$. Because $(f,f')$ is a contravariant pair
we have $X_f = X_{f'} \neq 0$ and because all cocircuits of a mutation are neighbours they have all the same $f'$-value.
\end{proof}

\begin{prop}\label{prop:NewMutationContravariantPair}
Let $\mathcal{O}$ be a Euclidean oriented matroid program of rank $\rk$ and let $(f,f')$ be a contravariant pair of elements (not being loops or coloops) of $\mathcal{O}$ with $f'$ lying in general position. Then the cocircuit graph $G^-_f$  of the program $(\mathcal{O},f',f)$ is empty.
If $\mathcal{O}$ has more than $\rk$ elements, the element $f$ has an adjacent mutation with all cocircuits of the mutation having $f' = +$ and the cocircuit $X$
of the mutation with $X_f \neq 0$ has $X_f = +$. The mutation remains a mutation in $\mathcal{O} \setminus f'$
\end{prop}

\begin{proof}
Because there is an element in general position in $\mathcal{O}$ not being a coloop, $\mathcal{O}$ has more then $\rk$ elements, hence 
Proposition \ref{prop:GfNotEmpty} yields that the cocircuit graph $G_f$ of the program $(\mathcal{O},f',f)$ is not empty.
A cocircuit having $f \neq 0$ and having $f' = +$ must have also $f = +$ hence the graph $G^-_f$ is empty which yields
$G^+_f$ being not empty.
The assumptions of Theorem \ref{theo:simplicialTopeMandel} are fulfilled and we get the desired mutation.
All cocircuits in that mutation have $f' = +$ hence Proposition \ref{prop:MutationsForContravariatPairs} yields that the mutation remains a mutation in
$\mathcal{O} \setminus f'$.
\end{proof}

\section{Las Vergnas oriented matroids}\label{section:LasVergnasOrientedMatroids}

We collect some results about Las Vergnas oriented matroids.
First, from Remark \ref{rem:dualOfR20MandelOrNot} we can conclude:

\begin{theorem}\label{theo:LasVergnasDualNotClosed}
The class of Las-Vergnas oriented matroids is not closed under duality.
\end{theorem}

If we look at the construction of the $R(20)$ in \cite{Richter-Gebert}, Chapter 2, we see that the oriented matroid after the second mutation process has an element with one
adjacent mutation. We have 

\begin{lemma}
There exist Las-Vergnas oriented matroids having an element with only one adjacent mutation.
\end{lemma}

We can apply Lemma  \ref{lemma:LexExtensionYieldsOneMutation} and Lemma \ref{lemma:LexExtensionleavesNonAdjacentMutations} from the Section before to show the following two Theorems.

\begin{theorem}
The class of Las-Vergnas oriented matroids is strictly bigger than the class of Mandel-oriented matroids.
\end{theorem}

\begin{proof}
We take the $R(20)$ oriented matroid again. It has one mutation-free element $f$.
Let $\mathcal{O}' = R(20)[f^-,e_1,e_2,e_3]$ a lexicographic extension of $R(20)$ with the elements $e_1,e_2,e_3$ arbitrarily chosen.
 Then because of Lemma \ref{lemma:LexExtensionYieldsOneMutation} the elements $f,f'$ are adjacent to a mutation in $\mathcal{O}'$.
 All the other mutations of $R(20)$ stay in $\mathcal{O}'$ because of Lemma \ref{lemma:LexExtensionleavesNonAdjacentMutations}.
Hence $\mathcal{O}'$ has no mutation-free element, it is a Las-Vergnas oriented matroid,
but it is not Mandel, because it has a minor that is not Mandel and because the class of uniform Mandel-oriented matroids is minor-closed, which is Theorem
\ref{theo:MandelOMsAreMinorClosed}. 
\end{proof}

We proved also the following:

\begin{theorem}
The class of Las-Vergnas oriented matroids is not minor-closed.
\end{theorem}

\section{Euclideaness and mutations}\label{section:EuclideanessAndMutations}

This section is an application of the first main result of our paper \cite{LexExt}. 
We cite it in two versions. First, we cite here Corollary 3.1 of \cite{LexExt}.  

\begin{corollary}\label{cor:lexExtStaysEucl2}
Let $\mathcal{O}$ be an oriented matroid of rank $\rk$ with groundset $E$ and let $e \in E$ being in general position. 
Let ${I} = [e_1, \hdots, e_{k}]$ with $k \le \rk$ be an ordered set of independent elements of $E \setminus e$ such that
$(\mathcal{O},e_1,e), (\mathcal{O} / \{e_1\},e_2,e), \hdots$, $(\mathcal{O} / \{e_1, \hdots, e_{k-1}\},e_k,e)$ are Euclidean oriented matroid programs.
Let $\mathcal{O}' = \mathcal{O} \cup p$ be the lexicographic extension $\mathcal{O}[e^+_1, \hdots, e^+_k]$.  
Then $(\mathcal{O}',p,e)$ is a Euclidean oriented matroid program.
\end{corollary}

\begin{remark}\label{rem:Rank4LexExtenRemainsEucl}
Be aware that all rang $3$ oriented matroids are Euclidean. Hence in the rank $4$ case the assumptions of the corollary reduce to $(\mathcal{O},e_1,e)$ being Euclidean.
 \end{remark}
 
 If $f$ is not in general position we need more assumptions.
 
\begin{theorem}[\cite{LexExt}, Theorem 3.1]\label{theo:lexExtStaysEucl2}
Let $\mathcal{O}$ be a Euclidean oriented matroid with $f \in E$ and let $\mathcal{O}' = \mathcal{O} \cup p$ be a lexicographic extension of $\mathcal{O}$.
Then $(\mathcal{O}',p,f)$ is a Euclidean oriented matroid program.  
\end{theorem}
  Together with the results of Section \ref{section:LexExtAndMutations} we can
   increase the estimation of the minimum number 
of adjacent mutations in Euclidean oriented matroids. 

\begin{theorem}\label{theo:3AdjacentMutations}
Let $\mathcal{O}$ be a simple connected Euclidean oriented matroid without having coloops of rank $\rk \ge 2$ with $n \ge \rk +2$ elements. 
In the case of rank $\le 3$ or corank $\le 3$ each element has at least $\rk$ adjacent mutations and $\mathcal{O}$ 
has at least $n$ mutations.
In the case of rank $> 3$ (with corank $> 3$) each element of $\mathcal{O}$ has 
at least $3$ adjacent mutations and  $\mathcal{O}$ 
has at least $n *3 / \rk$ mutations. 
\end{theorem}

\begin{proof} 
Oriented matroids in rank $2$ are realizable, also in corank $2$ because duals of realizable
oriented matroids stay realizable. There is nothing to show.
Hence, we assume rank and corank of $\mathcal{O}$ being $\ge 3$ hence we have $n \ge \rk +3$ elements. 
We show that each element $f$ (being not a coloop) has at least $3$ adjacent mutations (the rank $3$ case was already shown in in \cite{Levi26}, see \cite{OrientedMatroids}, Theorem 6.5.2 (i)).
We choose an element $g \neq f \in E$ with $(f,g)$ not being an inseparable pair. 
If there is no such an element, Lemma \ref{lem:ElementWithOnlyInseparableElIsRealizable} yields $\mathcal{O}$
being a realizable oriented matroid and we are done.
Like in the proof of Theorem \ref{theo:MandelIsLAsVergans} we use a lexicographic extension $\mathcal{O}' = \mathcal{O} \cup g' = \mathcal{O}[g^+,f^+,\hdots]$.
Then $g'$ is in general position and $(g',f)$ is also a separable pair.
Hence there are cocircuits $X,Y$ having $X_f = - Y_f \neq 0$ and $X_{g'} = Y_{g'} \neq 0$.
Then Theorem \ref{theo:lexExtStaysEucl2} yields that $(\mathcal{O},g',f)$ and $(\mathcal{O},g',e)$ for all $e \in E \setminus f$ are Euclidean oriented matroid programs  and the graphs $G_f^+$ and $G_f^-$ are not empty. 
Because of Theorem \ref{theo:simplicialTopeMandel}, $f$ has two adjacent mutations $M_1 = (f,e_1,e_2,e_3, \hdots)$ and $M_2 = (f,e_4,e_5,e_6, \hdots)$ (or simplicial topes $T^1$ and $T^2$). 
All vertices of $T^1$ have $g' = +$ and $f \in \{0,+\}$ and of $T^2$ have $g' = +$ and $f \in \{0,-\}$. The topes $T^1$ and $T^2$ remain different mutations in $\mathcal{O}$.
We assume (using reorientation) $T^1$ being a positive tope. Let $X$ be the cocircuit of $T^1$ with $X_f = +$ and let $Y$ be the cocircuit of $T^2$ with $Y_f = -$.

First, we assume having an element $e \neq f$ separating the topes $T^1$ and $T^2$, hence  the cocircuits of $T^1$ have $e \in \{0,+\}$ and the cocircuits of $T^2$ have $e \in \{0,-\}$. Because $e$ is not parallel to $f$ we have in each $T^1, T^2$ one cocircuit with $f = 0$ and $e \neq 0$.  
We take the lexicographic extension $\mathcal{O}' = \mathcal{O} \cup f' =  \mathcal{O}[f^+,e^-, \hdots]$. We have 
\begin{itemize}
\item a cocircuit $X^1$ in $T^1$ with $X^1_f = 0$ and $X^1_e = +$ having $X_{f'} = -$, 
\item the cocircuit $X$ in $T^1$ with $X_f = +$ having $X_{f'} = +$, 
\item a cocircuit $Y^2$ in $T^2$ with $Y^2_f = 0$ and $Y^2_e = -$ having $Y^2_{f'} = +$ and 
\item the cocircuit $Y$ in $T^2$ with $Y_f = -$ having $Y_{f'} = -$.
\end{itemize}

Second, if we do not have such an element $e$ the two mutations lie opposite to each other. The two mutations must then share their $f=0$ cocircuits. Otherwise in the contraction $\mathcal{O} / f$ we have again two mutations $T'^1$ and $T'^2$ with $T'^1$ being a positive tope there.
If we have a cocircuit $X$ being not in $T'^1$ it has $X_e = -$ for at least one element $e$ which is the separating element then.

Let then $X^1$ be a cocircuit of $T^1$ (hence also of $T^2$) having $X^1_f = X^1_{e_2} = 0$ and $X^1_{e_3} = +$
and let $X^2$ be a cocircuit of $T^1$ (hence also of $T^2$) having $X^2_f = 0$ and $X^2_{e_2} = +$.
We take the lexicographic extension $\mathcal{O}' = \mathcal{O} \cup f' =  \mathcal{O}[f^+,e_2^+,e_3^-, \hdots]$.
Then we have $X^1_{f'} = -$ and $X^2_{f'} = +$. In all cases $\mathcal{O}'$ remains Euclidean because of Theorem \ref{theo:lexExtStaysEucl2}.

Because $(f,f')$ are a contravariant pair, Proposition \ref{prop:NewMutationContravariantPair} yields that 
the cocircuit graph $G^+_f$ of $(\mathcal{O}',f',f)$ is not empty and it yields a mutation $M_3$ adjacent to $f$ with all cocircuits having $f' = +$.
But in $M_1$ as well as in $M_2$ there are cocircuits having opposite $f'$-values, hence the cocircuits from $M_3$ must differ from
these of $M_1$ or $M_2$.
The mutation $M_3$ stays a mutation in $\mathcal{O}$ that must be different from $M_1$ and $M_2$.
That means each element $f$ is adjacent to at least three mutations. If we count all elements, each mutation
is counted $\rk$ times hence $\mathcal{O}$ has at least $n * 3 / \rk$ mutations. In the rank $3$ case we have $n$ mutations
and also in the corank $3$ case because an oriented matroid has the same number of mutations than its dual.
%
\end{proof}

\begin{remark}
Before in the corank $3$ case it was only known that there are at least $n$ mutations but not that each element has at least $3$ adjacent mutations.
\end{remark}

We were not able to prove more for Euclidean oriented matroids. Also for the total number of mutations in a disconnected Euclidean 
oriented matroid, we can not get a such a nice result as in Lemma \ref{lem:disconnectedOMNumberMutations}.
Even for showing that there are $4$ adjacent mutations in the rank $4$ case we can not use the idea of our proof 
because we can not destroy in general three mutations with a lexicographic extension.
 To obtain $\rk$ mutations analogously to Shannons proof, see \cite{RoudneffSturmfels}, first proof of (1.1.2), or \cite{BokowskiMutationProblem}, Theorem 1.1,
 we would need to assume Levi's Intersection property for oriented matroids, see \cite{OrientedMatroids}, Definition 7.5.2, which is a stronger property than Euclideaness. 

\begin{definition}
An oriented matroid $\mathcal{O}$ of rank $\rk$ has {\em Levi's intersection property $IP_1$} if for every collection of $\rk -1$
hyperplanes $H_1, \hdots, H_{\rk -1}$ there is a single element extension $\tilde{\mathcal{O}} = \mathcal{O} \cup p$ such that $p$ is not 
a loop and $H_1 \cup p, \hdots, H_{\rk-1} \cup p$ are hyperplanes of $\tilde{\mathcal{O}}$.
\end{definition}

Moreover, we have to assume that the desired extension can be chosen to remain Euclidean. 
Then we can use this new extension as a separating element and apply
our proof idea of the lexicographic extension again to get one mutation more than before and using induction.
We could only use this idea in the uniform rank $4$ case.
We need a preliminary Proposition concerning mutations in the uniform case, before. 

\begin{prop}\label{prop:UniformOrientedMatroidsMutationsShareOnlyOneCocircuit}
In a uniform oriented matroid of rank $\rk$, two mutations share at most one cocircuit and a cocircuit is adjacent to at most two mutations.
If two mutations $M_1, M_2$ share a cocircuit $X$ they differ only in one element, hence $M_1 = [e_1, \hdots, e_{\rk-1}, f]$ and $M_2 = [e_1, \hdots, e_{\rk-1}, g]$
and we have $X^i_{e_i} = -Y^i_{e_i}$ for the base cocircuits for $i < \rk$ and $X^i_e = Y^i_e$ for all other elements $e \notin \{e_i,f,g\}$.
\end{prop}

\begin{proof}
First, let $X$ be a cocircuit with $z(X) = \{e_1, \hdots, e_{\rk-1}\}$
adjacent to a mutation $M_1 = [e_1, \hdots, e_{\rk-1}, f]$, we assume $M_1$ being a positive tope.
Let $M_2$ be a second mutation. If it shares the cocircuit $X$ with $M_1$, Proposition \ref{prop:ElementInGeneralPositionInvolvedInMutation} yields that all elements of $z(X)$
are adjacent to both mutations, hence we obtain
$M_2 = [e_1, \hdots, e_{\rk-1}, g]$.
 It is clear that the two mutations share only the
base cocircuit $X$. That proves the first statement.
Let $Y$ be the base cocircuit with $z(Y) = \{e_2, \hdots, e_{\rk - 1}, g\}$.
It must hold $Y_e = -$ for an element $e$ otherwise $Y$ would be adjacent to $M_1$ because it would be a positive cocircuit 
which is not possible because $Y_g = 0$ does not hold for any cocircuit of $M_1$.
We assume $Y_e = -$ for an $e \neq e_1$. But $Y$ is connected to $X$ via the edge $F$ with $z(F) = \{e_2, \hdots, e_{\rk-1}\}$. Let $X^1$ be the neighbour of $X$ in $M_1$ connected through that edge.
We have $X^1_e = +$ and elimination of $e$ between $X^1$ and $Y$ must yield $X$ hence we would have $X_e = 0$ which is impossible.
Hence we must have $Y_e = +$ for all elements $e \neq e_1$ and if must hold $Y_{e_1} = -$. 
This holds also for all other base cocircuits of $M_2$ which is the last statement of the Proposition. 
Now we assume there is a third mutation $M_3 =  [e_1, \hdots, e_{\rk-1}, h]$ adjacent to $X$. Let $Z$ be the base cocircuit having
$z(Z) = \{e_2, \hdots, e_{\rk-1}, h\}$. We have $Z \neq Y$, $Z_{e_1} = -$ and $Z_e = +$ for all other $e \neq e_1$. But $Y$ and $Z$ lie both on $F$
and we have $Y_g = 0, Y_h = +$ and $Z_h = 0, Z_g = +$ and $Y_e = Z_e$ for all other $e$. They are both different neighbours of $X$ hence it must hold $Z = X^1$
which is impossible.
\end{proof}

We come to our Theorem.

\begin{theorem}\label{theo:UniformIP2FourAdjacentMutations}
Let $\mathcal{O}$ be a uniform oriented matroid of rank $\rk = 4$ with $n > \rk + 3 = 7$ elements having the intersection property $IP_1$ and 
every extension obtained by that property can be chosen that it preserves Euclideaness. 
Then each element of $\mathcal{O}$ has at least $4$ adjacent mutations and  $\mathcal{O}$ 
has at least $n$ mutations. 
\end{theorem}

\begin{proof}
$\mathcal{O}$ is Euclidean because the $IP_1$-property implies Euclideaness.
We choose an arbitrary element $f$. From the Theorem before we have $3$ mutations adjacent to $f$. We show that there must be one adjacent mutation more.
Let $M_1, M_2, M_3$ be the three mutations with $X,Y,Z$ being the cocircuits having $f \neq 0$. 
The cocircuits $X,Y,Z$ must be different, we assume the contrary, $X = Y$ would yield $M_1 = [e_1, e_2, e_3, f]$ and $M_2 = [e_1, e_2, e_3, g]$
hence in $M_2$ we would not have cocircuits with $f = 0$. There is some case-checking to do.

First, we exclude the case that all three mutations have two elements $e,f$ in common.
Let $U^{i,h}$ be the $h$-base-cocircuit of mutation $M_i$ and let all $M_i$ have $e=+$. 
We have to check two cases. 
First, if all mutations have $f=\alpha$, we use the lexicographic extension $\mathcal{O} \cup f' = \mathcal{O}[f^+,e^{-\alpha}, \hdots]$.
We obtain $X_{f'} = Y_{f'} = Z_{f'} = \alpha$ and $U^{1,e}_{f'} = U^{2,e}_{f'} = U^{3,e}_{f'} = - \alpha$. 
Second, if two of the three mutations have $f=\alpha$, let $g$ be an element of the third mutation (we assume $M_3$).
Let $U^{3,g}_g = \beta$ hold.
We use the lexicographic extension $\mathcal{O} \cup f' = \mathcal{O}[f^+,e^{-\alpha}, g^{\alpha * \beta}, \hdots]$.
We obtain $X_{f'} = Y_{f'} = \alpha, Z_{f'} = -\alpha$ and $U^{1,e}_{f'} = U^{2,e}_{f'} = -\alpha$ and $U^{3,g}_{f'} = \alpha$. 
Hence in all mutations in the two cases, we have cocircuits with opposite $f'$-values.
The program $(\mathcal{O},f',f)$ stays Euclidean we proceed like in the proof before and obtain a fourth mutation.



Now we use the $IP_1$ and get an extension $\mathcal{O} \cup e$ through the cocircuits $X,Y$ and $Z$. 
We assume $\mathcal{O} \cup e$ remains Euclidean.
This extension can go through more than these three cocircuits of the mutations.
First, we assume it goes through more or equal than $5$ cocircuits.
That means more than two edges $F$ and $G$ with different zerosets of the mutations have $e=0$.
If we had $z(F) \cap z(G) = \emptyset$ in the original oriented matroid,
then it would hold $\rk(z(F) \cup z(G)) = 4$, hence intersecting those two edges would contradict submodularity. Hence we assume $\rk(z(F) \cup z(G)) = 3$
and $F$ and $G$ must have an element $g$ in common already in $\mathcal{O}$ (because $\mathcal{O}$ is uniform). 
The extension $e$ would then be parallel to $g$ which means
that we would have already $X_g = Y_g = Z_g$ in $\mathcal{O}$.
We have again the case that the three mutations have $2$ elements in common.

Second, we assume the extension goes only through the $3$ cocircuits $X,Y,Z$ of the mutations, we have $X_e = Y_e = Z_e = 0$. 
Then we can perturb these three cocircuits and
we get an element $e'$ such that all three mutations have cocircuits having $e' = +$ and $e' = -$.
The perturbation preserves Euclideaness of the programs $(\mathcal{O}, e', f)$. We show that later in Lemma \ref{lem:perturbationPreservesEuclideaness}.
We take a lexicographic extension $\mathcal{O} \cup e'' = \mathcal{O}[e'^+, \hdots]$ to bring the element in general position, 
the program $(\mathcal{O} \cup e'', e'', f)$ remains Euclidean because of Remark \ref{rem:Rank4LexExtenRemainsEucl}.
The contraction $\mathcal{O} / f$ is of rank $3$ hence Euclidean, therefore all the assumptions of Theorem \ref{theo:simplicialTopeMandel} are fulfilled.
We get a fourth mutation having $e'' = +$ for all of its cocircuits hence it must be different from the other three.

Third, the extension goes though $4$ cocircuits, hence $X,Y,Z$ and $W$.
There are two cases: the first one is that two mutations (say $M_1$ and $M_2$) are hanging together sharing the cocircuit $W$. They do not share 
a cocircuit $W'$ with $M_3$ because otherwise the extension would go also through $W'$, hence through $5$ cocircuits.
We assume $M_2$ being a positive tope hence $Y_f = +$ and $X_f = -$ and the cocircuit $Z$ of $M_3$ has $f = \alpha \neq 0$. 
The extension $\mathcal{O} \cup e$ through the cocircuits $X,Y$ and $Z$ goes also through the cocircuit $W$.
From the mutation $M_3$ we assume that only the cocircuit $Z$ has $Z_e = 0$, otherwise we obtain two edges in the extension, a case that we handled before.
Then we can perturb it which preserves Euclideaness of the program $(\mathcal{O} \cup e, e, f)$. This is again Lemma \ref{lem:perturbationPreservesEuclideaness}. 
The extension $e$ cuts the mutation $M_3$ then and we have $X_e = Y_e = 0$
and the two cocircuits $U^{1,e}_e = -$ and $U^{2,e}_e = +$ because of Proposition \ref{prop:UniformOrientedMatroidsMutationsShareOnlyOneCocircuit}.
We take the lexicographic extension $\mathcal{O} \cup f' = \mathcal{O}[e^+, f^-, \hdots]$. 
We obtain $Z_{f'} = Z_e$ and $U^{3,h}_{f'} = - Z_e$ for the other base cocircuits of $M^3$ because $e$ cuts the mutation $M_3$
and $X_{f'} = -X_f = +$ and $U^{1,e}_{f'} = U^{1,e}_e = -$ and $Y_{f'} = -Y_f = -$ and $U^{2,e}_{f'} = U^{2,e}_e = +$.
We obtain cocircuits with different $f'$-values and proceed like before. 

In the last case the mutations are not hanging together and 
we have one mutation, we assume $M_1$, having two cocircuits $X$ with $X_f = \alpha$ and $W = U^{1,e_1}$ with $X_e = U^{1,e_1}_e = 0$.
We have also $U^{1,e_1}_f = 0$ and let $U^{1,e_1}_{e_1} = \beta$. We perturb again $e$ to $e'$ at the other two mutations but leave it as it is at the mutation $M_1$.
Then we take the lexicographic extension $\mathcal{O} \cup e'' = \mathcal{O}[e'^+,f^-,e_1^{\alpha * \beta}, \hdots]$.
We have then $U^{1,e_1}_{e''} = \alpha$ and $X_{e''} = - X_f = - \alpha$. Again all three mutations have cocircuits with opposite $e''$-values.
We proceed like before. 
\end{proof}

We conjecture that this Theorem holds for all ranks and in general. To prove the general case,
it would be helpful if the following conjecture holds, which we think is also interesting in its own. 

\begin{conjecture}
If there is an extension $\mathcal{O} \cup p$ corresponding to a modular cut $\mathcal{M}$ of the corresponding matroid, and if 
$\mathcal{M}' \subseteq \mathcal{M}$ is also a modular cut, 
then there exist an extension $\mathcal{O} \cup p'$ corresponding to $\mathcal{M}'$.
\end{conjecture}

If the conjecture holds, we can always choose the extension yielded by the $IP_1$ such that it goes only through the three cocircuits $X,Y,Z$
which would avoid the case checking from above.
In the matroid case for matroids of rank $4$
Levi's Intersection property and Euclideaness are equivalent (see \cite{OrientedMatroids}, Proposition 7.5.6 (for matroids, not oriented matroids, it is wrong cited, see also  \cite{BachemWanka}, Proposition 4), but in the case of oriented matroids, 
the proof for that equivalence remains open. It is clear that we can use Euclideaness (property $IP_3$) to prove that we can always intersect two disjoint hyperplanes
with a point (or intersect two hyperplanes having a point together with a line). 
But because we do not know if this extension remains 
Euclidean (we give a counterexample in Section \ref{section:PreservationOfEuclideaness}) we cannot even intersect two disjoint hyperplanes with a line. Hence we cannot derive the $IP_2$ from $IP_3$ in the rank $4$ case in general.
We have only:

\begin{lemma} In the rank $4$ case
the $IP_3$ yields the $IP_2$ if the $IP_3$ is preserved in the extension and
the $IP_3$ yields the $IP_1$ if the $IP_3$ is preserved in the extension applied two times.
\end{lemma}

Finally we obtain 

\begin{theorem}
Let $\mathcal{O}$ be a uniform Euclidean oriented matroid of rank $4$ with $n > 7$ elements 
such that every parallel extension can be chosen to preserve Euclideaness three times. 
Then each element of $\mathcal{O}$ has at least $4$ adjacent mutations and  $\mathcal{O}$ 
has at least $n$ mutations. 
\end{theorem}

On the other hand all small classified Euclidean oriented matroids of rank $4$ (realizable or not) have at least $4$ mutations per element.
The question if there are Euclidean oriented matroids having an element with less than $\rk$ adjacent mutations stays open.

\section{Non-Euclidean Oriented Matroids}\label{section:Non-EuclideanOrientedMatroids}

In this Section we show some facts about directed cycles and very strong components in Non-Euclidean oriented matroids.

\begin{lemma}
A directed cycle in an oriented matroid program $(\mathcal{O},g,f)$ of rank $n$ has never edges on only one simplicial tope.
\end{lemma}

\begin{proof}
We assume that the conjecture holds in rank $n - 1$ and use induction.
Because all rank $\leq 3$ oriented matroids are Euclidean the lemma holds in that case trivially.
Let $X^1,X^2, \hdots, X^n$ be the vertices (cocircuits) of a rank $n$ simplicial tope and we assume we have a directed cycle
in $(\mathcal{O},g,f)$. All cocircuits in the cycle have $g = +$ and $f \neq 0$.
The cycle must contain all $n$ vertices otherwise we would have a directed cycle in a contraction violating the
induction assumption. 
W.l.o.g. we say 
\[ X^1_f = X^2_f = \hdots = X^n_f = - \text{ and let the cycle be } (X^1,X^2, \hdots, X^n). \]
Also let $e_i$ be the element where $X^1_{e_i} =  \hdots =  X^{i-1}_{e_i} = X^{i+1}_{e_i} = \hdots X^n_{e_i} = 0$. We assume (using reorientation) to have 
$X^i_{e_i} = +$ for all $i$.
 Since the cycle is directed we have  (let $n+1 = 1$): 
\[ El(-X^i,X^{i+1},g) = X^{i,i+1} \text{ with } (X^{i,i+1})_g = 0 \text{ and } (X^{i,i+1})_f = +.   \]
We obtain $X^{i,i+1}_{e_j} = 0$ for all $j \notin \{ i, i+1\}$
and $X^{i,i+1}_{e_i} = -$ and $X^{i,i+1}_{e_{i+1}} = +$. Then 
\[ El(X^{i+1},X^{i,i+1},f) = W^i \text{ with } W^i_g = + \text{ and } W^i_f = 0.\]
We obtain $W^i_{e_j} = 0$ for all $j \notin \{i,i+1\}$, $W^i_{e_i} = -$ and $W^i_{e_{i+1}} = +$. 
Then 
\[El(W^i,W^{i+1}, e_{i+1} = Z^i) \text{ with } Z^i_f = 0, Z^i_g = +, Z^i_{e_j} = 0 \text{ for all } j \notin \{i,i+2\}.\]
We have $Z^i_{e_ i} = -$ and $Z^i_{e_{i+2}} = +$. Then
\[ El(Z^i,W^{i+2}, e_{i+2}) = Z'^i \text{ with } Z'^i_f = 0, Z'^i_{e_j} = 0 \text{ for all }j \notin \{i,i+3\}.\]
We have $Z'^i_{e_ i} = -$ and $Z'^i_{e_{i+3}} = +$. 
We can go on like that and obtain a cocircuit $Z''^ i$ with 
$Z''^i_f = 0, Z''^i_g = +, Z''^i_{e_j} = 0$ for all $j \notin \{i,i-1\}, Z''^i_{e_ i} = -, Z''^i_{e_{i-1}} = +$. But
since $f \cup e_1 \cup \hdots \cup e_{i-2} \cup e_{i+1} \cup \hdots \cup e_n$ is a hyperplane in
the underlying matroid, we must have $Z''^i = \pm W^{i-1}$ contradicting $Z''^i_g = W^{i-1}_g = +$ and $Z''^i_{e_i} = -W^{i-1}_{e_i} = -$.  
\end{proof}

Unlike in the case of a simplicial tope, there exist directed cycles of oriented matroids programs having cocircuits only 
in one tope, see the $EFM(8)$ in \cite{OrientedMatroids}, Figure 10.4.7, but we can at least say that not all cocircuits of the tope can be involved.

\begin{lemma}
If the cocircuits of directed cycle in a uniform oriented matroid program $(\mathcal{O},g,f)$ are all contained in tope $\mathcal{T}$, 
then not all cocircuits of $\mathcal{T}$ can be in the cycle.
\end{lemma}

\begin{proof}
We assume the contrary.
By reorientation we assume $\mathcal{T}$ to be the feasible region of $(\mathcal{O},g,f)$. Because the cycle is directed 
none of its cocircuits has $g = 0$ or $f = 0$ hence the program is bounded. It is feasible because the cycle is not empty.
Hence Theorem \cite{OrientedMatroids}, 10.1.13 yields an optimal solution, which is a cocircuit in $\mathcal{T}$ that can not be part of the directed cycle which is impossible.
\end{proof}

We give a characterization of directed cycles in oriented matroids.

\begin{prop}
 From a directed cycle in $\mathcal{O}$, we can always derive 
another directed cycle $P = (X_1, \hdots, X_n = X_1)$ such that $X_i \circ X_{j}$ is not an edge
for all $j \neq i+1$. 
\end{prop}
 \begin{proof}
We assume there is such an edge $X^i \circ X^j$ with $i+1<j$. If the edge would have opposite direction, we obtain a smaller directed cycle $X_i, \hdots, X_j, X_i$.
Otherwise we go directly from $X_i$ to $X_j$ and get a smaller directed cycle. We do this as long as we can. 
\end{proof}

The next Lemma characterizes directed cycles in minimal non-Euclidean oriented matroids, 
see also the results in \cite{ExtensionSpaces}, after Lemma 4.10.
We call an oriented matroid {\em minimal non-Euclidean}, iff it has no non-Euclidean minor.

\begin{lemma}\label{lemma:characterizationMinimalNonEuclidean}
 Let $(\mathcal{O},g,f)$ be a uniform minimal non-Euclidean oriented matroid program and $P$ a directed cycle in the program.
 Then there is a directed closed path $P'$ containing the cocircuits of $P$ in $(\mathcal{O},g)$ such that the following holds:
 
\begin{enumerate}
 \item There is no $e$ such that for all $P_i \in P'$ holds $P_e = 0$.
 \item  There is no $e$ such that for all $P_i \in P'$ holds only $P_e = +$ or $P_e = -$.
 \item There is no $e$ such that for all $P_i \in P'$ holds only $P_e \in \{+,0\}$ or $P_e \in \{-,0\}$.
\item If there is an edge $P_i \circ P_{i+1}$ in $P'$ with $e \in z(P_i \circ P_{i+1})$ then $P'$ has two cocircuits $X$ and $Y$
 with $X_e = -Y_e \neq 0$. 
 \item Every $e$ is on an edge of $P'$.
 \end{enumerate}
\end{lemma}

\begin{proof}
The contradiction of the first two statements would violate the fact that $(\mathcal{O},g,f)$ is minimal non-Euclidean.
We could go to the contraction or to the deletion otherwise.
 Let there be an $e \in z(P^ i \circ P^{i+1})$ violating (iv). We assume $P^{i-1}_e = -$ for the predecessor of $P^i$ in the cycle.
  The path $P$ can be normalized to a path $P'$ with respect to $e$. (See here \cite{LexExt}, Definition 4.2 with $e=p$, and Lemma 4.5,
  because $\mathcal{O}$ is uniform all cocircuits with $e=0$ are new, see \cite{LexExt}, Theorem 2.2, the others old 
  and there are no undirected edges in the path.)
  Because $P'$ contains all cocircuits of $P$ it contains also the edge  $P^{i-1} \circ P^i$ which is directed from $P^{i-1}$ to $P^i$. 
 In the next edge $P^j \circ P^{j+1}$ in $P'$ where the cycle leaves the element $e$, the cocircuit $P^{j+1}$ must have $P^{j-1}_e = +$ then because $P'$
 is normalized.
 We do the normalization process as long as we can. This yields statement (iv). 
 Now we assume that (i),(ii) and (iv) hold already for the now normalized closed path $P$ and that 
 for all cocircuits of $P$ holds $e \in \{+,0\}$. Then we look at an edge $F = X \circ Y$ in $P$ with $X_e = +$ and $Y_e = 0$.
 We assume the edge is directed from $X$ to $Y$.
  The next neighbour $Z$ of $Y$ in $P$ can not have $Z_e = 0$ because that would violate Statement (iv).
  Also, $Z$ can not lie on $F$ because it would yield $Z_e = -$.
  Hence we need another edge $G$ different from $F$ on where the cycle is going on.
  We have $z(X) = z(F) \cup z(G) \cup e$ which is impossible because $\mathcal{O}$ is uniform. This yields iten (iii).
  We show the last item. Let again $(X,Y)$ be an edge in $P'$ with $X_e = -$ and $Y_e = 0$.
  If the next edge in the cycle is different from $F = X \circ Y$ it contains $e$ in its zeroset because $\mathcal{O}$ is uniform.
  Hence if $e$ is not lying on an edge in $P$ the next cocircuit in $P$ must be the other neighbour $Z$ of $Y$ lying in $F$.
  It holds $Z_e = +$. But then if we delete $e$ we have still the cocircuits $X,Z$ lying on $F$ being then a directed edge in $P \setminus e$.
We can omit $e$ and loose maybe that cocircuit but maintain the non-Euclideaness, having the edge $Y \circ Z$ conformal then.
\end{proof}

Here we recall the definitions of strongly connected and very strongly connected components, see \cite{ExtensionSpaces}, definitions before Example 3.5.

\begin{definition}
 We call two cocircuits of an oriented matroid program $(\mathcal{O},g,f)$ {\em very strongly connected} if they are connected by a directed closed path in both ways. This relation is an equivalence relation dividing the cocircuits in {\em components}. 
\end{definition}

The next Propositions are immediate.

\begin{prop}
 Every cocircuit of a non-isolated component is contained in a directed cycle.
\end{prop}

From Lemma \ref{lemma:characterizationMinimalNonEuclidean} we obtain a Theorem.

\begin{theorem}
 Let $(\mathcal{O},g,f)$ be a uniform minimal non-Euclidean oriented matroid program.
 Then if $P$ is a very strong component in the cocircuit graph of $(\mathcal{O},g)$ the following holds:
 For every $e \notin \{g,f\}$ of the groundset there is an edge in $P$ having $e$ in its zeroset and there 
 are cocircuits $X,Y$ in $P$ with $X_e = -Y_e \neq 0$.
\end{theorem}

We conjecture here that minimal Non-Euclidean oriented matroids have only one very-strong component in their cocircuit graphs, 
but we were not able to prove it.

\section{Preservation of Euclideaness}\label{section:PreservationOfEuclideaness}

We give a short survey of what is known about preservation of Euclideaness.
We know that it is preserved by taking minors (see \cite{OrientedMatroids}, chapter 10.5), taking the dual (see \cite{OrientedMatroids}, chapter 10.5), by exchanging target function and infinity (see \cite{LexExt},  Theorem 2.1) and the parallel extensions yielded by a Euclidean oriented matroid program preserve also Euclideaness in that program (see \cite{OrientedMatroids}, chapter 10.5). 
But these extensions can make other oriented matroid programs non-Euclidean. We show an example here.

\begin{lemma}
Let $\mathcal{O}$ be the oriented Vamos -Matroid. Let $l_1, l_2, l_3, l_4$ be the four lines in $\mathcal{M}(\mathcal{O})$ violating the bundle-condition, which means
all pairs of lines are coplanar except $l_3$ and $l_4$. Let $l_1 = \{a,b\}$ and $l_4 = \{c,d\}$. Let $\mathcal{O}' = \mathcal{O} \cup f$ be a lexicographic extension 
of $\mathcal{O}$ putting an additional point in general position on the line $\{b,d\}$.
Then $\mathcal{O}' \setminus d$ is Euclidean and has a parallel extension which makes the extended oriented matroid Non-Euclidean. 
\end{lemma}

\begin{proof}
Recall if we have a cocircuit $X$ in an oriented matroid program $(\mathcal{O}'',g,h)$ with $X_g = +$ and $X_h \neq 0$ and a parallel extension $f'$
going through $X$, then we have $f' \in \cl(z(X))$ and $f' \in \cl(g,h)$ hence we intersect the line $\cl(g,h)$ and the plane $z(X)$ with a point in the underlying matroid.
Now, $\mathcal{O} \setminus d$ is Euclidean and $\mathcal{O} \cup f \setminus d$ remains Euclidean because $f$ is a lexicographic extension here, too (lying in $l_1 \cup c$).
Then if we cut the line $\{b,f\}$ with the plane $l_2 \cup c$, we obtain $\mathcal{O}'$ which is Non-Euclidean. 
\end{proof}

In \cite{LexExt} we proved that the lexicographic extension of a matroid remains Euclidean.
We conjecture that Euclideaness is preserved for the connected sum of two oriented matroids described in \cite{Richter-Gebert}, Chapter 1. 
We will show this in a forthcoming paper in a special case for the oriented matroid $R(20)$ also given in \cite{Richter-Gebert}, Chapter 2. 
The reason to assume that conjecture is because we are dealing with inseparable elements (like in the lexicographic extension).
In the uniform case we can show a Lemma here, see here also \cite{ExtensionSpaces}, the proof of Proposition 4.4.

\begin{lemma}\label{lemma:EuclideanessRemainsForInseparablePairs}
Let $\mathcal{O}$ be a uniform oriented matroid and let $(f,f')$ be an inseparable pair of elements in $E$. 
Then $(\mathcal{O},g,f)$ is a Euclidean oriented matroid program iff $(\mathcal{O},g,f')$ is.
\end{lemma}

\begin{proof}
We assume $f$ and $f'$ being contravariant. We assume $(\mathcal{O},g,f')$ being Euclidean and
let $\mathcal{C}$ be a directed cycle in $(\mathcal{O},g,f)$.
All cocircuits of the cycle have $g = +$ and the same $f$-value, we assume $f = +$.
Let $(X,Y)$ be an edge in $\mathcal{C}$. If $X_{f'} = Y_{f'} \neq 0$ then $El(-X,Y,g) = Z$
with $Z_{f'} \neq 0$ and $Z_f \neq 0$ because $\mathcal{O}$ is uniform.
It holds $Z_{f'} = Z_f$ and the edge has the same direction in $(\mathcal{O},g,f')$.
Hence $\mathcal{C}$ would be a directed cycle in $(\mathcal{O},g,f')$ if all of its cocircuits had $f' \neq 0$
which is impossible. Because $f,f'$ are contravariant, Proposition \ref{prop:InseparablePairsAndNieghbours}  yields for each cocircuit $X$ having $X_{f'} = 0$
a neighbour $X'$ lying on the edge $F$ with $z(F) = z(X) \setminus f'$ having $X'_f = 0$ and $X'_{f'} = - X_f$.

Let now $(X,Y)$ be an edge in $\mathcal{C}$ with $X_{f'} \neq 0$ and $Y_{f'} = 0$. We have $X_f = Y_f = +$.
The other neighbour $Z$ of $Y$ in the edge $F = X \circ Y$ 
has then $Z_f = 0$ and $Z_g = +$ because $\mathcal{O}$ is uniform. Then $El(-Y,Z,g) = W$ with $W_f = -Y_f = -$ and $W_{f'} = Z_{f'} = -X_{f'}$. Then 
$El(-X,Y,g) = W'$ with $W'_{f'} = -X_{f'}$. Because $W = \pm W'$ must hold, we obtain $W = W'$ hence
the edge $(X,Y)$ is directed from $Y$ to $X$. In the cycle we can not come back from $X$ to another cocircuit having $f' = 0$.
Hence $\mathcal{C}$ contains only cocircuits having $f' = 0$.

Now, let $(X,Y)$ be an edge in $\mathcal{C}$ with $F = X \circ Y$ and $F_{f'} = 0$.
Let $X',Y'$ be the corresponding neighbours of $X$,$Y$ on the edges $z(X) \setminus f'$ and $z(Y) \setminus f'$ that have $X'_f = Y'_f = 0$.
It can not be that $X' = Y'$ because then $X' \setminus f$ would be also a cocircuit contradicting $\mathcal{O}$ being uniform.
We have $z(X' \circ Y') = z(X') \cap z(Y') = (z(X) \cap z(Y)) \cup f \setminus f'= F \cup f \setminus f' $, hence $X',Y'$ are comodular.
Proposition \ref{prop:InseparablePairsAndNieghbours} yields $X'_{f'} = -X_f = -Y_f = Y'_{f'}$ hence $f' \notin sep(X',Y')$.
If there would be another element $s \in sep(X',Y')$ then $s \in z(X) \setminus z(Y)$ would contradict that $\mathcal{O}$ is uniform,
 $s \in z(X \circ Y)$ would yield $s \in z(X' \circ Y')$ which is impossible, hence it must hold $X_s = -Y_s \neq 0$ contradicting that $X,Y$ are conformal.
 We obtain that $X',Y'$ are conformal, hence neighbours.
  
Then $El(-X,Y,g) = Z$ and $El(-X',Y',g) = Z'$. It is clear that $Z$ and $Z'$ must be neighbours again,
 hence we have $Z_{f'} = 0$ and $Z'_f = 0$ and Proposition \ref{prop:InseparablePairsAndNieghbours} yields again $Z'_{f'} = - Z_f$. 
  The edge $(X',Y')$ has the opposite direction in $(\mathcal{O},g,f')$ as $(X,Y)$ in $(\mathcal{O},g,f)$.
 We obtain a directed cycle in $(\mathcal{O},g,f')$ which is impossible.
 \end{proof}

%
%

Euclideaness remains if we take the direct sum of two oriented matroids.

\begin{lemma}
The direct sum of two Euclidean oriented matroids remains Euclidean.
\end{lemma}

\begin{proof}
Let $\mathcal{O}^1$ and  $\mathcal{O}^2$ be two Euclidean oriented matroids and
let $\mathcal{O}$ be their direct sum. 
We check cases for the program $(\mathcal{O},g,f)$.
First, if $g \in \mathcal{O}^1(E)$ and $f \in \mathcal{O}^2(E)$,
a cocircuit having $g = +$ must have all elements of $\mathcal{O}^2(E)$ in its zero set and vice versa.
Hence the program has only cocircuits with $f = 0$ in the graph $(\mathcal{O},g)$ hence is Euclidean.
Second if we have $f,g$ coming both from $\mathcal{O}^1(E)$ the cocircuits of $(\mathcal{O},g)$
are the same like the cocircuits in the program $(\mathcal{O}^1(E),g,f)$ which is Euclidean. 
\end{proof}

As a corollary we obtain:

\begin{corollary}
The union of two Euclidean oriented matroids remains Euclidean.
\end{corollary}

\begin{proof}
The construction of the union of two oriented matroids uses only the operation of direct sum, lexicographic extension and deletion, see \cite{OrientedMatroids}, the paragraph before Theorem 7.6.2 and the paragraph after Example 7.6.5. All these operations preserve Euclideaness.
\end{proof}

In general Euclideaness is not preserved by mutation-flips in oriented matroid programs, see a counterexample
in \cite{ExtensionSpaces}, Example 3.5 and the following but if either the target-function or the infinity (not both) are adjacent to the mutation, it is.
Recall that for uniform oriented matroids mutation-flips correspond to simplicial topes, see \cite{OrientedMatroids}, Theorem 7.3.9.
We need a proposition before.

\begin{prop}\label{prop:CocircuitOfMutationNotInDirectedCycle}
 Let $(\mathcal{O},g,f)$ be a uniform oriented matroid program of rank $\rk$ with $g$ being not a loop and $f$ being not a coloop.
 Let $M = [e_1, \hdots, e_r]$ be a mutation in $\mathcal{O}$ with $M \cap \{f,g\} \neq \emptyset$  
 then no cocircuit of the mutation lies in a directed cycle in $(\mathcal{O},g,f)$.
\end{prop}

\begin{proof}
All cocircuits of a directed cycle have $g = +$ and $f \neq 0$.
If $f$ and $g$ are both in $M$ then every cocircuit of $M$ has either $g=0$ or $f=0$ or both hence the proposition holds in that case trivially.
Now let $f \in M$ ($f = e_{\rk}$) and $g \notin M$. Then all cocircuits of $M$ have $f = 0$ except one cocircuit $X$ having $X_{e_1} = X_{e_2} = X_{e_{\rk-1}} = 0$.
We assume w.l.o.g. that $X_f = +$. Then it is easy to see that $X$ has no predecessor in $(\mathcal{O},g,f)$ having $f = +$, all edges going to other cocircuits
$Y$ with $Y_f = +$ are directed from $X$ to $Y$ hence $X$ can not be involved in a directed cycle.
\end{proof}

\begin{lemma}\label{lemma:MutationFlipOnlyTargetFunctionInvolvedStaysEuclidean}
 Let $(\mathcal{O},g,f)$ be a Euclidean uniform oriented matroid program  of rank $\rk$ with $g$ being not a loop and $f$ being not a coloop.
 Let $M = (f,e_2, \hdots,e_{\rk})$ be a mutation in $\mathcal{O}$ with $g \notin M$.
 Let $\mathcal{O}'$ be an oriented matroid derived by a mutation flip of $M$.
 Then $(\mathcal{O}',g,f)$ remains Euclidean.
\end{lemma}

\begin{proof}
Let $(X,Y)$ be an edge of the cocircuit graph of the program with $X,Y$ not lying in $M$.
Then $El(X,Y,g) = Z$. Because $Z_g = 0$ the cocircuit $Z$ is not adjacent to the mutation $M$ hence it does not change its values in $\mathcal{O}'$, 
see \cite{OrientedMatroids}, Corollary 7.3.6 and Definition 7.3.4, the set $\mathcal{U}$ is in our case the set of cocircuits adjacent to $M$. 
The direction of the edge $(X,Y)$ stays the same in $(\mathcal{O}',g,f)$. 
We assume a directed cycle in $(\mathcal{O},g,f)$. Because of Proposition \ref{prop:CocircuitOfMutationNotInDirectedCycle}
no cocircuits of the mutation are involved in the cycle. Hence all directions of the edges stay the same in  $(\mathcal{O}',g,f)$
and we have a directed cycle there.
\end{proof}

We receive a Corollary from that. We call the {\em mutation-graph} of a class of oriented matroids the graph
having the oriented matroids as vertices and edges between two of them if they are connected via a mutation-flip. 
We call an oriented matroid {\em totally Non-Euclidean} if it has no Euclidean oriented matroid programs.

\begin{corollary}\label{cor:TotallyNonEuclideanOMsMoreThenOneMutationFlipFromLinearOnes}
 Let $\mathcal{O}$ be a totally Non-Euclidean uniform oriented matroid of rank $4$.
 Then if it is connected in the mutation-graph to a Euclidean uniform oriented matroid it
 is connected by at least three mutation-flips. 
 \end{corollary}
 
 \begin{proof}
 Let $\mathcal{O}^1$ be the linear uniform oriented matroid. Let $M_1 = (1,2,3,4)$ be the first mutation.
 Then after the flip we obtain $\mathcal{O}^1$ with maybe the programs $(\mathcal{O}^1,g,f)$ be non-Euclidean with
 $\{g,f\} \subset M_1$ or  $\{g,f\} \cap M_1 = \emptyset$. The programs with $f \in M_1$ and $g \notin M_1$ (or vice versa) stay Euclidean
 because of Lemma \ref{lemma:MutationFlipOnlyTargetFunctionInvolvedStaysEuclidean}.
 Now we take the next mutation $M_2$. 
 Then the programs with $f \in M_2 \setminus M_1$ and $g$ in $M_1 \setminus M_2$ (and vice versa) stay Euclidean
 as well as the programs with $f \in M_1 \cap M_2$ and $g \in M^c_1 \cap M^c_2$. 
 Either the first or the second case or both must appear. Hence, we have still Euclidean oriented matroid programs.
We will need a third mutation-flip to make these non-Euclidean.  
With three mutation-flips we can theoretically reach totally Non-Euclideaness e.g. if all three mutations disjoint.
\end{proof}
 
That corollary yields that all uniform rank $4$ OMs with maximal $8$ points are not totally Euclidean. (Of course that follows already directly from the classification of all these oriented matroids, see \cite{bokowski1990classification} with a little computerhelp.)
Moreover, we characterize all $8$-point uniform Non-Euclidean rank $4$ oriented matroids as counterexamples for the other cases
of mutation-flips where Euclideaness is not preserved.

\begin{lemma}
Let $\mathcal{O}$ be a uniform $8$-point non-Euclidean oriented matroid of rank $4$. Then it has a linear mutant via a flip of a mutation $M$
and non-Euclidean oriented matroid programs $(\mathcal{O},a,b)$ where $a,b \in M$ and $(\mathcal{O},c,d)$ with $c,d \notin M$.
\end{lemma}

\begin{proof}
All these oriented matroids have a linear mutant, see  \cite{bokowski1990classification}. Because of Lemma \ref{lemma:MutationFlipOnlyTargetFunctionInvolvedStaysEuclidean} we have $(\mathcal{O},a,b)$ non-Euclidean with $a,b \in M$ or $a,b \notin M$. But we obtain always both cases because in all these oriented matroids all elements are involved in at least one non-Euclidean program (this is yielded by computerhelp). 
\end{proof}

Perturbation of oriented matroids can also preserve Euclideaness, we show a special case here
which we needed before in Theorem \ref{theo:UniformIP2FourAdjacentMutations}.

\begin{lemma}\label{lem:perturbationPreservesEuclideaness}
Let $\mathcal{O}$ be a uniform oriented matroid of rank $4$ having a mutation $M = [f,e_2,e_3,e_4]$.
Let $X$ be the $f$-base cocircuit. Let $\mathcal{O} \cup e$ be an extension of $\mathcal{O}$ such that $X_e = 0$
and all other base-cocircuits of $M$ have $e = +$. Let all programs $(\mathcal{O}, e, g)$ for all $g \notin \{e_2,e_3,e_4\}$ of the groundset of $\mathcal{O}$
be Euclidean oriented matroid programs. Then it is possible to perturb $e$ to get an oriented matroid $\mathcal{O}$ such that $X_e = -$ and the $e$-value of all other cocircuits stay the same. The programs $(\mathcal{O}', e', g)$ remain Euclidean.
\end{lemma}

\begin{proof}
First, we have to show that the perturbation is always possible. This is \cite{OrientedMatroids}, Theorem 7.3.1.
Second, we show that Euclideaness remains.
Let $g$ be an element of the groundset, not part of the mutation, 
Because $M = [e',e_2,e_3,e_4]$ is a mutation in the program $(\mathcal{O}', e', g)$,
Proposition \ref{prop:CocircuitOfMutationNotInDirectedCycle} yields that the cocircuit $-X$ is never used in a directed cycle in that program.
If we have another edge $F = Y \circ Z$ in $(\mathcal{O}, e', g)$, elimination of $e'$ between $-Y$ and $Z$
yields a cocircuit $W'$ with $W'_{e'} = 0$ and $W'_g = \alpha$. 
Elimination of $e$ between $-Y$ and $Z$ in $(\mathcal{O}, e, g)$
yields a cocircuit $W$ with $W_e = 0$ and $W_g = \beta$.
Now if $W \neq X$ we have $W' = W$ and the direction of $F$ stays the same in both programs.
If $W = X$ in $(\mathcal{O}, e, g)$ we have now $W'$ being a cocircuit of the new mutation.
But all cocircuits of that mutation have the same $g$-value, the direction of the edge remains also in that case.
\end{proof}

\section{Non-Euclidean but Mandel oriented matroids}\label{section:NonEuclideanButMandel}

Here we show that Inclusion (3) of Equation \ref{equ:inclusions} is strict, there are many Mandel oriented matroids that are non-Euclidean. 
We know that the lexicographic extenion preserves Euclideaness, 
 but here we are dealing with non-Euclidean oriented matroids, 
and have to look carefully at the assumptions of the corresponding theorems in  \cite{LexExt}.   
We will need Corollary \ref{cor:lexExtStaysEucl2} from Section \ref{section:EuclideanessAndMutations} again
and we formulate a lemma for a special case with weaker assumptions for uniform oriented matroids.

\begin{lemma}\label{lem:lexExtStaysEucl3}
Let $\mathcal{O}$ be a uniform oriented matroid of rank $\rk$ with groundset $E$. 
Let ${I} = [f,e_2, \hdots, e_{\rk}]$ be an ordered set of elements of $E$ such that
$(\mathcal{O},f,e_2)$, $(\mathcal{O} / \{e_2\},f,e_3), \hdots$, $(\mathcal{O} / \{e_2, \hdots, e_{\rk-1}\},f,e_{\rk})$ are Euclidean oriented matroid programs.
Let $\mathcal{O}' = \mathcal{O} \cup f'$ be the lexicographic extension $\mathcal{O}[f^+,e^+_2, \hdots, e^+_{\rk}]$.  
Then $(\mathcal{O}',f,f')$ is a Euclidean oriented matroid program.
\end{lemma}

\begin{proof}
Let $I = [e_2, \hdots, e_{\rk}]$ and let $\mathcal{O}'' = \mathcal{O} \cup e'$ be the lexicographic extension $\mathcal{O}[e^+_2, \hdots, e^+_{\rk}]$. 
We show that $(\mathcal{O}',f,f')$ is Euclidean iff $(\mathcal{O}'',f,e')$ is Euclidean which proves everything together with Corollary \ref{cor:lexExtStaysEucl2}.
Let $(X,Y)$ be an edge of  a directed cycle $\mathcal{C}$ in $(\mathcal{O}',f,f')$.
Then $X_f = Y_f = +$ and $X_{f'} = Y_{f'} \neq 0$ and $El(X,Y,f) = Z$ with $Z_f = 0$.
It holds $Z_{f'} \neq 0$ because $\mathcal{O}'$ is uniform, 
hence we have $Z_{f'} = Z_{e_i} \neq 0$ where $i$ is the first index such that $Z_{e_i} \neq 0$.
But this yields also $Z_{e'} = Z_{e_i}$ in $(\mathcal{O}'',f,e')$.
The edge has the same direction in $(\mathcal{O}',f,f')$ as in $(\mathcal{O}'',f,e')$.
That shows everything.
\end{proof}

 The next proposition is obvious.

\begin{prop}\label{prop:contractionMutationFlipStays}
Let $\mathcal{O}$ be a uniform oriented matroid and $\mathcal{O}_M$ derived from $\mathcal{O}$ by a flip of a mutation $M$.
Let $g$ be an element not in $M$. Then $\mathcal{O}/g = \mathcal{O}_M/g$
\end{prop}

\begin{proof}
Because $g$ is not part in the mutation and because $\mathcal{O}$ is uniform all cocircuits in $M$ have the same $g$-value, which is $\neq 0$.
These cocircuits do not exist anymore in the contractions.
\end{proof}

\begin{prop}\label{prop:MutAndLexExtExchange}
Let $\mathcal{O}$ be a uniform oriented matroid of rank $r$ with groundset $E$ and with $f,g \in E$. 
Let $M = [f, e_2, \hdots, e_r]$ be a mutation and all elements of $M_1$ directed versus the mutation. 
Let all cocircuits of the mutation have $g = +$.
\begin{itemize}
\item Let  $\mathcal{O}_{f'} = \mathcal{O} \cup f' = \mathcal{O}[f^+, g^-, e^{\alpha_3}_3, \hdots, e_r^{\alpha_r}]$ be a lexicographic extension \\ of $\mathcal{O}$.
\item We have a mutation $M' = [f', e_2, \hdots, e_r]$ in $\mathcal{O}_{f'}$. 
\item Let $\mathcal{O}_{f', M'}$ be the oriented matroid obtained from $\mathcal{O}_{f'}$ by a mutation-flip of $M'$.
Then $M$ and $M'$ are mutations in $\mathcal{O}_{f', M'}$.
\item Let  $\mathcal{O}_M$ be the oriented matroid obtained from $\mathcal{O}$ by a mutation-flip of $M$. 
\item Let $\mathcal{O}_{M,f'} = \mathcal{O}_M \cup f' = \mathcal{O}_M[f^+, g^+, e^{- \alpha_3}_3, \hdots, e_r^{- \alpha_r}]$
be a lexicographic extension of $\mathcal{O}_M$. 
\item Then $M'$ is a mutation in $\mathcal{O}_{M,f'}$.
\item Let $\mathcal{O}_{M,f',M'}$ be the oriented matroid obtained from $\mathcal{O}_{M,f'}$ by a mutation-flip of $M'$.
Then $M$ and $M'$ are mutations in $\mathcal{O}_{M,f',M'}$.
\end{itemize}
Then $\mathcal{O}_{f',M'}$ and $\mathcal{O}_{M,f',M'}$ are isomorphic under the isomorphism which maps $f$ to $f'$ and vice versa 
and is the identity for all other elements of $E$. The proposition holds also if the $e_3, \hdots, e_r$ for the lexicographic extensions are arbitrarily chosen in $E \setminus \{f,g\}$.
\end{prop}

 \begin{proof} 
  Lemma \ref{lemma:MutationsAndLexExtensions}
  yields that $M'$ is a mutation in $\mathcal{O}_{f'}$ and in $\mathcal{O}_{M,f'}$.
  Let  $E_M = \{e_2, \hdots, e_r\}$. 
  We have $M = [f,E_M]$, and $M' = [f', E_M]$.
 The proof works analogously to the proof of Lemma \ref{lemma:LexExtensionIsomorphism}. 
 The lexicographic extensions used here are like in the Lemma. For cocircuits outside of the mutations (hence not involved in the reorientation processus of the mutation-flips) the mapping from $f$ to $f'$ and vice versa is the same as
in Lemma \ref{lemma:LexExtensionIsomorphism}.
  Only the cocircuits in the mutations $M$ and $M'$ need to be considered.
  First, all these cocircuits share the same values for elements in $E \setminus E_M$ and they have all $g=+$.
  We consider the $f$-base cocircuit $X$ of the mutation $M$ in $\mathcal{O}$. It has $X_f = +$ in $\mathcal{O}$, hence
$X_{f'} = +$ in $\mathcal{O}_{f'}$ hence $X_f = +$ and $X_{f'} = -$ in $\mathcal{O}_{f',M'}$.
It has $X_f = -$ in $\mathcal{O}_M$ hence $X_f = X_{f'} = -$ in  $\mathcal{O}_{M,f'}$ and $X_f = - $ and  $X_{f'} = +$ in $\mathcal{O}_{M,f',M'}$.
Hence $X$ has in $\mathcal{O}_{f',M_2}$ and in $\mathcal{O}_{M,f',M'}$ the same values if we exchange $f$ and $f'$.

Now let $Y^i$ be an $e_i$-base-cocircuit of the mutation $M$ in $\mathcal{O}$. It has $Y^i_f = 0$ and $Y^i_{e_i} = +$ in $\mathcal{O}$ and 
$Y^i_{f'} = -$ in $\mathcal{O}_{f'}$. Hence we have
\[Y^i_f = 0, Y^i_{e_i} = + \text{ and } Y^i_{f'} = - \text{ in } \mathcal{O}_{f',M'}.\]
It has $Y'^i_f = 0$ and $Y'^i_{e_i} = -$ in $\mathcal{O}_M$ and $Y'^i_f = 0$, $Y'^i_{f'} = +$ and $Y'^i_{e_i} = -$ in $\mathcal{O}_{M,f'}$. Hence we have 
\[Y'^i_{f'} = +, Y'^i_{e_i} = - \text{ and } Y'^i_f = 0 \text{ in } \mathcal{O}_{M,f',M'}.\]

Now let $Z^i$ be an $e_i$-base-cocircuit of the mutation $M'$ in $\mathcal{O}_{f'}$. It has $Z^i_{f'} = 0$, $Z^i_f = +$ and $Z^i_{e_i} = +$ in $\mathcal{O}_{f'}$.
We have
\[ Z^i_f = +, Z^i_{e_i} = -  \text{ and } Z^i_{f'} = 0 \text{ in } \mathcal{O}_{f',M'}. \]

Now let $W^i$ be an $e_i$-base-cocircuit of the mutation $M'$ in $\mathcal{O}_{M,f'}$. It has $W^i_{f'} = 0$, $W^i_f = -$ and $W^i_{e_i} = -$ in $\mathcal{O}_{M,f'}$.
Hence it has 
\[ W^i_{f'} = 0, W^i_{e_i} = + \text{ and } W^i_f = - \text{ in } \mathcal{O}_{M,f',M'}. \]
Now if we exchange $f$ and $f'$ in $\mathcal{O}_{M,f',M'}$ We see that the cocircuits $Y$ and $W$ and the cocircuits $Y'$ and $Z$
have the same values. Because mutation-flips do not change the underlying matroid, we have in $\mathcal{O}_{f',M_2}$ and in $\mathcal{O}_{M,f',M'}$ 
the same number of cocircuits which proves all.
 \end{proof}

 \begin{lemma}\label{lem:ConstructionOfMandelExtension}
Let $\mathcal{O}$ be a uniform non-Euclidean oriented matroid of rank $\rk$
having a mutation $M = (f,e_2,\hdots,e_{\rk})$ with all elements of $M$ oriented versus the mutation. Let $\mathcal{O}/f$ be Euclidean and
let $\mathcal{O}_M$ be a Euclidean mutant derived from $\mathcal{O}$ via the mutation flip of the mutation $M$.
Let $g \notin M$ and we assume all cocircuits of $M$ have $g = +$. 
Let $\mathcal{O}_{f'} = \mathcal{O} \cup f'$ be the lexicographic extension $\mathcal{O}[f^+,g^-,e_3^-, \hdots,e_{\rk}^-]$. 
It has the mutation $M' = (f',e_2,\hdots,e_{\rk})$.
Let $\mathcal{O}_{f',M'}$ be the oriented matroid derived from $\mathcal{O}'$  via the mutation-flip using $M'$.
Then $(\mathcal{O}_{f',M'},e,f')$ is Euclidean for all $e$ in the groundset of $\mathcal{O}$ and we have $\mathcal{O}_{f',M'} \setminus f' = \mathcal{O}$.
 \end{lemma}
 
 \begin{proof}
 Lemma \ref{lemma:MutationsAndLexExtensions} yields that $M'$ is a mutation in $\mathcal{O}_{f'}$. Because of $g \notin M$, Proposition \ref{prop:contractionMutationFlipStays} yields that
 $\mathcal{O}/g = \mathcal{O}_M/g$ is Euclidean.  
Because $\mathcal{O}_M$ is Euclidean, the programs  $(\mathcal{O}_M,f,e)$ are Euclidean for all $e$ and  
Lemma \ref{lemma:MutationFlipOnlyTargetFunctionInvolvedStaysEuclidean} yields that $(\mathcal{O},f,e)$
with $e \notin M$ (hence also for $e = g$) is Euclidean. Because  $\mathcal{O}/f$ is Euclidean all assumptions of Corollary \ref{cor:lexExtStaysEucl2} (set $e_1 = f$) are fulfilled and we obtain that $(\mathcal{O}_{f'},f',e)$, hence $(\mathcal{O}_{f'},e,f')$ is Euclidean for all $e \notin M$. 
Also, because $(\mathcal{O},f,g)$ is Euclidean and because  $\mathcal{O}/g$ is Euclidean the assumptions from Lemma \ref{lem:lexExtStaysEucl3} (set $e_2 = g$) are fulfilled
and $(\mathcal{O}_{f'},f,f')$ is Euclidean.
Hence we have $(\mathcal{O}_{f'},e,f')$ Euclidean for all $e \notin M'$ and Lemma \ref{lemma:MutationFlipOnlyTargetFunctionInvolvedStaysEuclidean}
yields that $(\mathcal{O}_{f',M'},e,f')$ is Euclidean for all $e \notin M'$, hence for all $e \notin M$ and for $e = f$. 

Now let $e \in M \setminus f$. Because $\mathcal{O}_M$ is Euclidean the program $(\mathcal{O}_M,f,e)$ is Euclidean. 
We extend it lexicographic. Let $\mathcal{O}_{M,f'} = \mathcal{O}_M \cup f' = \mathcal{O}_M[f^+,g^+,e_3^+,,\hdots,e_{\rk}^+]$. 
The program $(\mathcal{O}_{M,f'}, f', e)$ remains Euclidean because of Corollary \ref{cor:lexExtStaysEucl2} (set $e_1 = f$), hence also $(\mathcal{O}_{M,f'}, e,f')$.
Then, because  $f$ and $f'$ are inseparable, Lemma \ref{lemma:EuclideanessRemainsForInseparablePairs} yields that $(\mathcal{O}_{M,f'}, e, f)$ remains Euclidean.
Again $M' = (f', e_2,\hdots, e_{\rk})$ is a mutation in $(\mathcal{O}_{M,f'}, e, f)$. We have $e \in M'$ and $f \notin M'$, hence Lemma   \ref{lemma:MutationFlipOnlyTargetFunctionInvolvedStaysEuclidean} yields that in the flipped oriented matroid $\mathcal{O}_{M,f',M'}$
the program $(\mathcal{O}_{M,f',M'},e,f)$ remains Euclidean. Proposition \ref{prop:MutAndLexExtExchange} yields that $\mathcal{O}_{M,f',M'}$ is 
isomorphic to $\mathcal{O}_{M,f'}$ if we exchange $f'$ and $f$. The last statement is obvious.
\end{proof}
  
It would be nice to get rid of the contraction-assumption of this lemma. This assumption holds always in rank $4$ case 
and for minimal non-Euclidean oriented matroids.
From that lemma follow directly two theorems.

\begin{theorem}
A uniform oriented matroid $\mathcal{O}$ of rank $4$ that is reachable by a Euclidean uniform oriented matroid via one mutation-flip is Mandel.
\end{theorem}

\begin{proof}
We may assume  $\mathcal{O}$  being non-Euclidean, hence we have more than $\rk + 3$ elements.
Because we are in rank $4$ all contractions of $\mathcal{O}$ are Euclidean.
All assumptions (using maybe reorientation) of Lemma \ref{lem:ConstructionOfMandelExtension} are fulfilled and we are done.
\end{proof}

The same proof works for the next theorem.

\begin{theorem}
A minimal non-Euclidean uniform oriented matroid that is reachable by a Euclidean uniform oriented matroid via one mutation-flip is Mandel.
\end{theorem}

We conclude that every uniform rank $4$ oriented matroid with $8$ points is Mandel since they are all reachable from linear oriented matroids by one mutation-flip. It shows that Inclusion (3) is proper. Because there is a uniform $8$ point oriented matroid of rank $4$ having an element with only 
$3$ adjacent mutations, we have an example of a rank $4$ Mandel oriented matroid with an element with only three adjacent mutations.

\section{concluding remarks}\label{section:concludingRemarks}
Let $L$ be the minimum number of adjacent mutations to an element of a class of matroids.
Summarizing the results of this paper we have for connected rank $\rk$ oriented matroids:

\begin{theorem}
\begin{enumerate}
\item $L( \frak{O} ) = 0.$
\item $L(\frak{O}_{LasVergnas}) = 1$.
\item $1 \le L(\frak{O}_{Mandel}) \le  3$ if $\rk \ge 3$.
\item $L(\frak{O}_{Euclidean}) =  \rk$ if $\rk \le 3$.
\item $3 \le L(\frak{O}_{Euclidean}) \le  \rk$ if $\rk > 3$.
\item $L(\frak{O}_{IP2',uniform}) = 4$ if $\rk = 4$.
\item $L(\frak{O}_{realizable}) =  \rk$.
\end{enumerate}
\end{theorem}

We summarize the open questions raised in this paper.

\begin{itemize}
\item Is the class of Mandel oriented matroid in generally minor-closed? 
\item \label{question:MandelDualClosed} Is the class of Mandel oriented matroid closed under duality? 
\item Are there Euclidean oriented matroids with less than $\rk$ mutations?
\item Are there Mandel oriented matroids with less than $3$ mutations?
\item Are there Mandel oriented matroids with more mutation-flip distance than $1$ to a Euclidean oriented matroid? 
\end{itemize}
For totally Non-Euclidean oriented matroids we have an example which is LasVergnas and one which is not.
One question stays open:
\begin{itemize}
\item \label{question:TotallyNonEuclMandel} Is a non Las-Vergnas oriented matroid always totally Non-Euclidean? 
On the other hand are there Mandel totally Non-Euclidean oriented matroids?
\end{itemize}
We would like to answer it for the two other known not LasVergnas-examples.
A more special question is:
\begin{itemize}
\item Is the dual of the $R(20)$ Mandel?
\end{itemize}
If the answer of that question is positive then the class of Mandel oriented matroids would be minor closed.
\begin{itemize}
\item Are the three intersection properties equivalent for oriented matroids in rank $4$?
\end{itemize}
This is shown in the matroid case and open in the oriented matroid case.
A positive answer would yield $L = 4$ for uniform Euclidean oriented matroids of rank $4$.
\begin{itemize}
\item Can we always intersect two hyperplanes in a rank $4$ oriented matroid?
\end{itemize}
This is a crucial question which would show the equivalence from before and is also interesting of its own.
\begin{itemize}
\item Are the three intersection properties not equivalent for oriented matroids in rank $\ge 5$?
\end{itemize}
Again this is shown for the matroid case. We do not know if the counterexamples given there are orientable.
We showed that the lexicographic extension creates at least one mutation. We ask:
\begin{itemize}
\item How many mutations creates a lexicographic extension?
\item Is Euclideaness preserved by the construction of the conneccted sum in \cite{Richter-Gebert}, Chapter 1?
\end{itemize}

\bibliography{MutationsAndNonEuclideaness}
\bibliographystyle{plain} 
\nocite{*}

\end{document}